\documentstyle [11pt,amstex,amscd]{article}
\def\bd{\partial}

\def\dim{{\rm dim\; }}

\def\text#1{{\em #1}}

\def\la{\lambda}

\def\be{\begin{equation}}
\def\ee{\end{equation}}
\def\bear{\begin{eqnarray}}
\def\eear{\end{eqnarray}}
\def\best{\begin{eqnarray*}}
\def\eest{\end{eqnarray*}}
\def\pf{{\bf Proof}: }

\newtheorem{theorem}{Theorem}[section]
\newtheorem{prop}[theorem]{Proposition}

\newtheorem{lemma}[theorem]{Lemma}
\newtheorem{cor}[theorem]{Corollary}

\newtheorem{defn}[theorem]{Definition}

\newtheorem{remark}[theorem]{Remark}
\newenvironment{rem}{\begin{remark}\rm}{\end{remark}} 

\newtheorem{example}[theorem]{Example}
\newenvironment{ex}{\begin{example}\rm}{\end{example}}

\def\non{\noindent}
\def\pf{\non {\bf Proof. }}
\def\qed{\nopagebreak \hskip .1in { $\Box$ }\penalty10000 %
\hskip\parfillskip \par  }

\def\ra{\rightarrow}

\def\r#1{\right#1}
\def\l#1{\left#1}

\def\ma#1{\mathop {#1} \limits}

\def\Si{\Sigma}

\def\ti{\times}

\def\P{{ \Bbb P}}
\def\Q{{ \Bbb Q}}
\def\cx{{ \Bbb C}}

\def\wt#1{\widetilde{#1}}
\def\wh#1{\widehat{#1}}
\def\ov#1{\overline{#1}}

\def\M{{\cal M}}

\def\longra{\longrightarrow}
\def\H1{{\cal H}^1}
\def\si{\sigma}
\def\Si{\Sigma}
\def\dim{\mbox{dim }}

\def\O{{\cal O}}
\def\A{{\cal A}}
\def\phi{\varphi}
\def\longra{\longrightarrow}
\def\z{{\cal Z}}
\def\x{{\cal X}}
\def\y{{\cal Y}}
\def\oz{\ov\z}
\def\oy{\ov\y}
\def\D{{\bf D}}
\def\L{{\cal L}}

\title{\bf Topological Recursive Relations in $H^{2g}(\M_{g,n})$}

\author{ Eleny-Nicoleta Ionel\thanks{partially supported by
the N.S.F. and by a Sloan research fellowship} \\ University of Wisconsin \\
 Madison, WI 53706 }

\date{}
\addtocounter{section}{0}
\begin{document}

\maketitle

\vskip.15in

\begin{abstract}
We show that any degree at least $g$ monomial in descendant or 
tautological classes  vanishes on $\M_{g,n}$ when $g\ge 2$. This 
generalizes a result of 
Looijenga and proves a version of Getzler's conjecture. The method we 
use is the study of the 
relative Gromov-Witten invariants of $\P^1$ relative to two points combined 
with the degeneration formulas of [IP1]. 
\end{abstract}



\vskip.4in

Let $\ov\M_{g,n}$ be the Deligne-Mumford compactification of the 
moduli space $\M_{g,n}$ of genus $g$ smooth curves with $n$ (distinct) 
marked points. Let $L_i\ra \ov\M_{g,n}$ be the relative cotangent
bundle at the marked point $x_i$; the fiber of $L_i$ over 
$(\Si,x_1,\dots,x_n)$ is the cotangent space to $\Si$ at $x_i$. The 
first Chern class of this bundle is denoted $\psi_i=c_1(L_i)$ and 
is sometimes called a (gravitational) descendant. If  
$\pi:\ov\M_{g,n+1}\ra \ov\M_{g,n}$ is the map that forgets the last 
marked point then $\kappa_a=\pi_*(\psi_{n+1}^{a+1})$ is called a  
tautological 
class (or Mumford-Morita-Miller class); since 
$\kappa_a\in H^{2a}(\ov\M_{g,n})$ we define its degree to be $a$,
while the degree of each $\psi_i$ equals 1. 
\smallskip

In \cite{l1} Looijenga proved that in the Chow group $\A^*({\cal C}^n_g)$ 
a product of descendant classes of degree at least $g+n-1$ vanishes, where 
${\cal C}^n_g$ is the moduli space of smooth genus $g$ curves with $n$ (not 
necessarily distinct) points. In particular, in $\M_{g,0}$ 
any degree $g-1$ monomial in tautological classes vanishes. However,
with the above definition of tautological classes, this not 
true anymore in $\M_{g,n}$, for $n\ge 1$ (for example in $\M_{2,1}$  
$\kappa_1=\psi_1\ne 0$).

In this paper, we obtain the following generalization of Looijenga's result:
\begin{theorem} When $g\ge 2$, any product of degree at least $g$ (or
at least $g-1$ when $n=0$) of descendant or tautological classes 
vanishes when restricted to $H^*(\M_{g,n},\Q)$. 
\label{T.main}
\end{theorem}
Note that when $g\le 1$, is has been known for a long time that
$\psi_j$ and $\kappa_a$ with $a\ge 1$ vanish on $\M_{g,n}$. 

\medskip

The proof of Theorem \ref{T.main} is a simple consequence of the degeneration 
formula for relative Gromov-Witten invariants (cf. \cite{ipa}). The
idea is to start with the moduli space $\y_{d,g,n}$ of degree 
$d$ holomorphic maps from a smooth genus $g$ surface with $n$ marked points
into $\P^1$ which have a fixed ramification pattern over $k$ marked
points in the target $\P^1$. In Section 1 we describe the structure of
$\y_{d,g,n}$ and that of its compactification $\oy_{d,g,n}$. The 
relatively stable map compactification $\oy_{d,g,n}$ is closely related to
both the space of admissible covers (introduced by Harris-Mumford
in  \cite{hmu}) and the space of twisted covers (recently defined  by 
Abramovich-Vistoli in \cite{av}). Moreover, it comes with
two natural maps $st$ and $q$ that record respectively the
domain and the target of the cover. One of the key ideas of the paper
is then to pull back by $q$ known relations in the cohomology of the target
and then push them forward by $st$ to get relations in the cohomology
of the domain. So we need to know that the space $\oy_{d,g,n}$ carries
a fundamental class (over $\Q$) and that it satisfies Poincare
duality. The discussion in Section 1 shows this assertion, so in particular  
$st_*\oy_{d,g,n}$ defines a cycle in $\ov\M_{g,n}$; the codimension of
this cycle is at most $g$ when $\oy_{d,g,n}$ is a 2-point ramification cycle
(i.e. all but two of the branch points are simple). 

We next choose the degree $d$ and a 2-point ramification cycle so that the 
stabilization map $st:\oy_{d,g,n}\ra \ov\M_{g,n}$ has finite, nonzero degree. 
Theorem \ref{T.1} then shows that any 
product of descendants on the domain is a linear combination of
(generalized) 2-point ramification cycles on $\ov\M_{g,n}$. There are  three
main ingredients in its proof. We first  relate the
relative cotangent bundle of the domain to the pull back via $q$ of
the relative cotangent bundle of the target. Next, it is known that
when  genus is zero then (nontrivial) products of  descendants are 
Poincare dual to boundary cycles $D$ in $\ov\M_{0,k}$ 
(see for example \cite{keel}). This relates a 
product of descendants on the domain  to cycles of type $st_*q^* D$, 
and the degeneration formula (\ref{deg.form})  completes the proof 
of Theorem \ref{T.1}. 

Corollary \ref{P.induction2} then implies that the Poincare dual of 
any degree $m$ 
product of descendant and tautological classes can be written as a linear
combination generalized 2-point ramification cycles of codimension
$m$. But the codimension of a 2-point ramification cycle is at most
$g$; Proposition \ref{P.B.is.bd} proves that the cycles of codimension exactly
$g$ vanish on $\M_{g,n}$, thus finishing the proof of Theorem
\ref{T.main}. All degenerations used in this paper are in fact linear 
equivalences, so an algebraic-geometric proof of the degeneration formula
(\ref{deg.form}) would in fact give not only the vanishing in
cohomology, but also in the Chow ring, as in Looijenga's Theorem. 

\medskip

From Theorem \ref{T.1} we see that the 2-point ramification cycles on
$\ov\M_{g,n}$ generate a subring that contains descendant and
tautological classes. In fact, we believe that this subring is not
larger then the one generated by descendant, tautological classes and
their pullbacks by the attaching maps of the boundary strata of $\ov\M_{g,n}$. 
At least when restricted to $\M_{g,n}$, the
arguments in Section 7 of \cite{mu} easily extend to show that any
2-point ramification constraint can be expressed as a polynomial in
descendants and tautological classes.  It's not clear at this moment
how to generalize this argument to the compactification $\ov
\M_{g,n}$.

On the other hand, when the genus is low ($g\le 5$) one can prove that 
all 2-point ramification constraints appearing in Theorem \ref{T.1} 
are in fact polynomials in only descendant and  tautological classes 
supported on the boundary. Moreover, the coefficients of this
polynomial can be determined by keeping track of the coefficients in 
(\ref{deg.form}). Relations expressing products of descendant classes
as polynomials in descendant and tautological classes 
supported on the boundary are known as topological recursive 
relations (TRR). The $g=0$ and $g=1$ TRR's were known classically. 
In genus 2, Mumford derived a  
formula for $\psi_1^2$  and  Getzler (\cite{g}) for  $\psi_1\psi_2$. 
In the same recent paper \cite{g}, Getzler made the conjecture that
for any genus $g$ there are degree $g$ TRR's. 

When the genus is 3 for example, Theorem \ref{T.main} implies the
following new relations (modulo boundary terms): 
$\psi_1^2\psi_2=\psi_1\psi_2\psi_3=0$ 
(as Getzler conjectured), plus the unexpected relation 
$\kappa_1\psi_1\psi_2=0$. Unfortunately, if we keep track of the boundary 
terms, the number of terms in the TRR increases very fast as the genus grows. 
The genus 0 and genus 1 TRR have 1 and 2 terms 
respectively, but the genus 2 TRR in \cite{g} has 18 boundary terms.
We leave  the actual TRR formulas in low genus ($3\le g\le 5$) for 
another paper.
\bigskip

Note that the degree $g$ is the lowest degree in which 
one could hope that some {\em monomial} in descendants would vanish 
on $\M_{g,n}$. 
The reason is that the class $\psi_2\dots \psi_n \la_g\la_{g-1}$ vanishes on 
$\bd \ov\M_{g,n}$ (where $\la_i=c_i(E)$ are the Chern classes 
of the Hodge bundle), while Faber's  conjecture (\cite{fa}), which also 
agrees with Virasoro predictions (see  \cite{gp}) gives
\best
\psi_1^{a_1}\psi_2^{a_2+1}\dots \psi_n^{a_n+1}\la_g\la_{g-1}= 
{(2g-3+n)!\over (2a_1-1)!!(2a_2+1)!!\dots (2a_n+1)!!}
\cdot {|B_{2g}|\over 2^{2g}g(2g-1)!}\ne 0
\eest
when $\ma\sum_{i=1}^n a_i= g-1$. On the other hand, for large genus,
there are most likely  lower degree (homogeneous) polynomials in 
descendants which vanish 
on $\M_{g,n}$.   
\bigskip

While this paper was under revision, the author heard a conjecture
made by Vakil \cite{v}. He essentially conjectured that in the Chow
group, any degree $m$ monomial in $\kappa$ and $\psi$ classes on
$\ov\M_{g,n}$ is pullback from the strata with at least $m+1-g$ genus
0 components. In the cohomology group, this conjecture follows 
immediately from the results of this paper, and was added 
as the final Proposition \ref{P.ravi}. As mentioned above, an
algebraic-geometrical proof of the degeneration formula 
(\ref{deg.form}) would also  give the result in the Chow group.  
\vskip.4in

\section{The space of relatively stable covers}
\bigskip

We start by defining a space of degree $d\ge 1$, Euler characteristic
$\chi$ covers of $\P^1$  with prescribed ramification pattern over several
points of $\P^1$. The ramification indices at each point $p\in \P^1$ 
will be encoded by an ordered sequence of positive 
multiplicities $I=(s_1,\dots, s_\ell)$. For any such $I$, we define
\best
\ell(I)=\ell\qquad \deg I=\ma\sum_{i=1}^\ell s_i \qquad
|I|=\ma\prod_{i=1}^\ell s_i.  
\eest
We also allow some of the points in the inverse image of $p$ to be  
marked points on the domain.  
\begin{defn} 
Consider $I_1,\dots,I_k$  ordered sequences of  
multiplicities with $\deg \;(I_j)=d\ge 1$ for all $j$, and 
let $N_1, \dots ,N_k$ be an ordered partition of the set 
$\{x_1,\dots, x_n\}$ (where some of the $N_j$'s might be
empty). For all $j=1,\dots,k$ assume that 
 $0\le \ell(N_j)\le  \ell(I_j)$, where
$\ell(N_j)$ denotes the cardinality of $N_j$. We define 
\bear
\Xi_{d,\chi}  \l(\prod_{j=1}^k  b_{I_j}(N_j)\r)
\label{hatX.0}
\eear
to be the infinite dimensional manifold consisting of data 
$(f,\Si, J,x_1,\dots,x_n;p_1,\dots,p_k)$ such that:
\begin{enumerate}
\item[(i)] $J$ is a complex structure on $\Si$, a smooth two dimensional
real manifold (not necessarily connected) with Euler characteristic $\chi$;
\item[(ii)]   $x_1,\dots x_n$ and 
$p_1,\dots, p_k$ are distinct points on $\Si$ and respectively $\P^1$; 
\item[(iii)]  $f:(\Si,J)\ra \P^1$ is a degree $d$ holomorphic map, which has
moreover positive degree on each component of $\Si$; 
\item[(iv)]  for each $j=1,\dots,k$, there exist distinct points
$(x_{n_{ij}}\,)_{i=\ell(N_j)+1,\dots,\ell(I_j)}$ on $\Si$, distinct from 
$x_1,\dots, x_n$ such that 
\best
f^{-1}(p_j)=
\ma\sum_{i=1}^{\ell(I_j)} s_{ij}\, x_{n_{ij}}
\eest
(i.e. $f$ is ramified at $x_{n_{ij}}$ of index $s_{ij}$), 
where $I_j=(s_{ij})_{i=1,\dots ,\ell(I_j)}$ and
$N_j=(x_{n_{ij}})_{i=1,\dots, \ell(N_j)}$. 
\end{enumerate}
By convention, the space (\ref{hatX.0}) is empty when 
$\ell(N_j)> \ell(I_j)$ or $\deg I_j\ne d$.
\end{defn}
We say that $b_{I_j}(N_j)$ describes the ramification pattern of $f$
over the point $p_j\in \P^1$. Note that when $\deg I_j>\ell(I_j)$ the
point  $p_j$ is a branch point of multiplicity $\deg I_j-\ell(I_j)$. 
For example 
$b_{2,1^{d-2}}(x_1)$ means that $x_1$ is a simple ramification point
while $b_{1^d}(x_1,x_2)$ means that $x_1$ and $x_2$ are conjugate points of 
the cover. 

In this context, we can think of
$b_{I_j}(N_j)$ as imposing a $(\deg I_j-\ell(I_j)+\ell(N_j))$-dimensional 
condition on a generic degree $d$ covering map 
$f:(\Si,x_1,\dots,x_n)\ra (\P^1,p_1,\dots,p_k)$. In particular, we
usually work with ramification patterns $b_{I_j}(N_j)$ that
satisfy $\deg I_j-\ell(I_j)+\ell(N_j)\ge 1$.

The space (\ref{hatX.0}) 
has several components, depending on the
topological type of the domain $\Si$; the component corresponding to a
fixed $\Si$ will be denoted by 
\best
\Xi_{d,\Si} \l(\prod_{j=1}^k  b_{I_j}(N_j)\r)
\eest
 
\begin{defn} 
The groups $Diff(\Si)$ of diffeomorphisms of $\Si$ and 
$Aut(\P^1)=PGL(2,\cx)$ of automorphisms of $\P^1$   
act on $\Xi_{d,\Si}\l( \ma\prod_{j=1}^k b_j(N_j)\r)$ by 
\best
(g,h)\cdot(f,\Si,J,x_1,\dots,x_n,p_1,\dots,p_k)= (h\circ f\circ
g,\Si, g^*J,g^{-1}(x_1),\dots, g^{-1}(x_n),h(p_1),\dots,h(p_k))
\eest
where $g\in Diff(\Si)$ and $h\in Aut(\P^1)$. Consider the two quotients
\best
\wh \x_{d,\Si}\l( \ma\prod_{j=1}^k b_{I_j}(N_j)\r)=
\l.\Xi_{d,\Si}\l(\ma \prod_{j=1}^k b_{I_j}(N_j)\r)
\r/{Diff(\Si)}
\eest 
and
\bear
\x_{d,\Si}\l( \ma\prod_{j=1}^k b_{I_j}(N_j)\r)=
\l.\Xi_{d,\Si}\l(\ma \prod_{j=1}^k b_{I_j}(N_j)\r)
\r/{Diff(\Si)\ti Aut(\P^1) }
\label{def.X}
\eear
The latter is called the moduli space of {\em smooth} degree $d$ covers of
$\P^1$ by $\Si$ with  ramification pattern $b_{I_j}(N_j)$ at points 
$p_j\in \P^1$
for $j=1,\dots,k$. The corresponding union of spaces $\x_{d,\Si}$ over 
different topological types $\Si$ with the same Euler characteristic 
$\chi$ is denoted 
by
\best
\x_{d,\chi}\l( \ma\prod_{j=1}^k b_{I_j}(N_j)\r).
\eest
\label{Def.X}
\end{defn}
\vskip-.2in 

An element $f\in \x_{d,\chi}$ is an equivalence class of triples
consisting of a smooth domain $C=(\Si,j,x_1,\dots,x_n)$, the (marked) 
target $(\P^1,p_1,\dots, p_k)$ and the covering map. The groups
$Diff(\Si)$ and $Aut (\P^1)$ have induced actions on the domain and
respectively the target. Therefore the space $\x_{d,\chi}$   comes with 
two natural projections
\begin{equation}
\begin{CD}
\ov\M_{0,k}@<{q}<<\x_{d,\chi} \l(\ma\prod_{j=1}^k  b_{I_j}(N_j)\r)@>{st}>>
\wt\M_{\chi,n}
\end{CD}
\label{def.diag}
\end{equation}
defined by $q(f)=(\P^1,p_1,\dots,p_k)$ and $st(f)=C$, where
$\wt\M_{\chi,n}$ is the moduli space of complex structures with $n$
marked points on a possibly disconnected
curve with Euler characteristic $\chi$. In fact, after choosing some 
ordering the $m$ components of $\Si$ we 
see that
\best
\wt \M_{\chi,n}=\ma\bigsqcup_{m=1}^\infty \;\;
\l(\; \ma\bigsqcup \;\ov\M_{g_1,n_1}\ti\dots\ti \ov\M_{g_m,n_m}\;\r)/S_m
\eest  
where the second union is over all  $g_i$, $n_i$ and distributions of
the $n$ marked points on the $m$ components such that 
$\ma\sum_{i=1}^m (2g_i-2)=\chi$, $\ma\sum_{i=1}^m n_i =n$; the 
symmetric group $S_m$ acts
by permuting the $m$ components. 

Restricting to a fiber of $q$ in the fibration
(\ref{def.diag}) gives us a corresponding moduli space of covers  
with prescribed ramification pattern at $k$ {\em fixed points } in
$\P^1$, denoted 
\best
\x_{d,\Si} \l(\ma\prod_{j=1}^k B_{I_j}(N_j)\r).
\eest
The $k$ points are suppressed in the notation for convenience. 
\medskip

\begin{rem} Since the degree of the covering map $f$ is required to be 
positive on each  component of $\Si$ and the group $Diff(\Si)$ acts on 
$\Xi_{d,\Si}$ with finite stabilizers then
 $\wh\x_{d,\Si}\l( \ma\prod_{j=1}^k b_{I_j}(N_j)\r)$ 
has a natural orbifold structure of dimension
\best
\dim \wh \x_{d,\Si} \l(\prod_{j=1}^k  b_{I_j}(N_j)\r)&=&2d-\chi(\Si)+k+n-
\sum_{j=1}^k(\deg(I_j)-\ell(I_j)+\ell(N_j))
\\
&=&2d-\chi(\Si)-\sum_{j=1}^k(\deg(I_j)-\ell(I_j))+k
\eest
When moreover $k\ge 3$ then $Aut(\P^1)$ also acts with finite 
stabilizers, so in this case the quotient  
$\x_{d,\Si}\l( \ma\prod_{j=1}^k b_{I_j}(N_j)\r)$  
is naturally an orbifold  of dimension 
\best
\dim  \x_{d,\Si} \l(\prod_{j=1}^k  b_{I_j}(N_j)\r) =2d-\chi(\Si)-
\sum_{j=1}^k(\deg(I_j)-\ell(I_j))+k-3
\eest
When $k\le 2$, $Aut(\P^1)$ has a $3-k$ dimensional subgroup which acts
trivially and so $\x_{d,\Si}$ still has an orbifold structure, but of
dimension $2d-\chi(\Si)-\ma\sum_{j=1}^k(\deg(I_j)-\ell(I_j))$. 

Similarly, when $2g-2+n\ge 1$, the moduli space $\ov\M_{g,n}$ has an orbifold
structure of dimension $3g-3+n$ (obtained by adding Pyrm structures
as described in \cite{l2}), while when $2g-2+n\le 1$
it has a (nonstandard) orbifold structure of dimension $g$. More
precisely, for $n\le 3$ we have $\ov\M_{0,n}=\ov\M_{0,3}=pt$ and
similarly $\ov\M_{1,0}=\ov\M_{1,1}$.    
\end{rem}
\medskip

The space $\x_{d,\chi}$ also comes with a collection of intrinsic line
bundles. Denote by $L_{x_i}\ra\wt\M_{\chi,n}$ and $L_{p_j}\ra \ov\M_{0,k+r}$ 
the relative cotangent bundles at the marked points $x_i$ and $p_j$
respectively. Next, let $\L_{x_i}\ra \x_{d,g}$ be the relative
cotangent bundle to the (unstabilized) domain $C$ at the marked point
$x_i$ and $\L_{p_j}=q^*L_{p_i}\ra \x_{d,g}$ be the relative 
cotangent bundle to the target $\P^1$ at $p_j$. The fiber at 
$f\in \x_{d,g}$ of $\L_{x_i}$ is $T^*_{x_i}C$ while that of 
$\L_{p_j}$ is $T^*_{p_j}\P^1$. To eliminate the possibility of
confusion, throughout this paper $x$ will denote a marked point of the 
domain and $p$ will denote a marked point of the target.

\bigskip

We next want to compactify $\x_{d,\chi}$ so that the maps in
 the diagram (\ref{def.diag}) extend continuously and so that
$st_*\ov\x_{d,\chi}$ defines a cycle in $\wt\M_{\chi,n}$.  For that, we use
the {\em relatively stable} maps compactification of the space of
smooth holomorphic maps into $\P^1$ relative to the collection of 
marked points
$\{p_1,\dots, p_k\}$ in the (target) $\P^1$ (cf. Section 6 of
\cite{ip4}). This compactification is similar in spirit to the usual
`stable maps into $\P^1$'  compactification (as described for example
 in \cite{p}) but it is much finer. The
difference is that not only the domain can bubble (or equivalently
gets rescaled) when for example two marked points start colliding, but 
also the target $\P^1$ gets rescaled around $p_j$ when a ghost component
(i.e. collapsed component) starts forming or the points $p_j$ get too close
to each other.

The strata in the usual stable map compactification that have ghost 
components not only have the wrong dimension (their obstruction bundle
comes from the Hodge bundle of the ghost domain), but more 
importantly, if the ghost component is sent to $p$, the ramification  
constraint above the point $p$ becomes undefined. Making the target 
 bubble yields in the limit a holomorphic map to a  degenerate $\P^1$, but
without any ghost components over $p$. 
\medskip

More precisely, consider a sequence $(f_n)$ of smooth 
degree $d$ stable holomorphic maps to $\P^1$ that have a fixed
ramification pattern $b_I(N)$ above $p$. Suppose that their usual
stable map limit $f$ has some ghost components $C_2$ over $p$. Let
$f_1:C_1\ra \P^1$ be the restriction of $f$ to the other components of
$C$ and let $b_S$ be its ramification pattern above $p=p_0$ (in
general $S\ne I$). After rescaling the target $\P^1$ around $p$ (and
passing to a subsequence) we obtain in the limit a second nontrivial
cover $f_2:C_2\ra \P^1$ that has the same ramification pattern $b_S$
over $p_\infty$, and fewer (if any) ghost components over $p$.  If
$f_2$ still has ghost components over $p$, we continue rescaling.
Otherwise, $f_2$ has the ramification pattern $b_I(N)$ over $p$ and
all together the limit map is a degree $d$ cover 
\bear 
f=f_1\cup
f_2:C_1\ma\cup_{y_i^1=y_i^2\atop i=1\dots \ell}C_2\ra \P^1\ma
\cup_{p_0=p_\infty} (\P^1,p) 
\eear
of a degenerate $\P^1$ (with an
ordinary double point). The cover $f$ has no ghost components over $p$
or the nodal point $p_0=p_\infty$, and $f_1^{-1}(p_0)=\sum s_i y^1_i$,
$f_2^{-1}(p_\infty)=\sum s_i y^2_i$ so $f_1$, $f_2$ have the same
ramification pattern $b_S$ over the node $p_0=p_\infty$.
\medskip

To have a good compactification of 
$\x_{d,\chi}\l(\ma\prod_{j=1}^k b_{I_j}(N_j)\r) $ we
must use the rescaling process around at least all the points $p_j$
for $j=1,\dots ,k$, so that the limit map still satisfies the 
ramification constraints $b_{I_j}(N_j)$ at the points $p_j$. However,
things become simpler to describe if there are no other branch
points. In what follows we restrict our attention to the moduli space
of stable maps where {\em all} the branch points are marked:
\begin{defn} Define a moduli space of possibly disconnected {\em
smooth} covers 
\bear
&&\z_{d,\chi}  \l(\prod_{j=1}^k  b_{I_j}(N_j) \r) 
\;\;\ma=^{def}\;\;
\x_{d,\chi}  \l(\prod_{j=1}^k  b_{I_j}(N_j) \; ( b_{2,1^{d-2}})^r
\r)
\quad 
\label{def.Z}
\eear
where the last $r$ branch points are simple and ordered, with $r$ given by
\bear
r=2d+\chi-\sum_{j=1}^k (\deg I_j-\ell(I_j))
\label{def.r}
\eear
When $\chi=2-2g$ let 
\bear
&&\y_{d,g}  \l(\prod_{j=1}^k  b_{I_j}(N_j) \r) \subset 
\z_{d,\chi}  \l(\prod_{j=1}^k  b_{I_j}(N_j) \r) 
\quad 
\label{def.Y}
\eear
denote the subspace of connected covers.
\label{D.mod}  
\end{defn}  
Recall from  Definition \ref{Def.X} that an element $f$ of the space 
$\x_{d,\chi}$ is a triple consisting of a
marked domain, marked target and a degree $d$ covering map with a 
specified ramification pattern at marked points in the target. 
All the images of marked points in the domain are marked; some of the 
preimages of the marked points of the target might also marked. 
However, there are possibly many unmarked ramified points mapping to 
marked or unmarked points of the target. An element $f$ of the 
space $\z_{d,\chi}$ has the extra property that all its branch points 
are marked in the target, and in particular the ramification pattern
of $f$ is completely determined. 
\bigskip

 Moreover, when $k+r\ge 3$ the space  
$\z_{d,\chi}\l(\ma\prod_{j=1}^k  b_{I_j}(N_j) \r) $ has a canonical
orbifold structure of dimension 
\bear
\dim \z_{d,\chi}  \l(\prod_{j=1}^k  b_{I_j}(N_j) \r) =2d-\chi-
\sum_{j=1}^k(\deg(I_j)-\ell(I_j))+k-3=r+k-3.
\label{dim.Z}
\eear
When $k+r\le 2$ Lemma \ref{0.dim} below  shows that the space 
$\z_{d,\chi}\l(\ma\prod_{j=1}^k b_{I_j}(N_j) \r) $ is  0
dimensional. 

\begin{lemma} Consider the space 
$\y=\y_{d,g}  \l(\ma\prod_{j=1}^k  b_{I_j}(N_j) \r)$ and let $r$ be as in 
Definition \ref{D.mod} while $n=\ma\sum_{j=1}^k \ell(N_j)$. If 
$2g+n\ge 3$ then
$k+r\ge 3$. Moreover, if $k+r\le 2$, then $\y$ consists of only one
element; the domain of this cover is an unstable $g=0$ curve and 
the covering is totally ramified at two points. 
\label{0.dim}
\end{lemma}
\pf When $k=2$ relation (\ref{def.r}) becomes
$r=2g-2+\ell(I_1)+\ell(I_2)$. So $r>0$ unless $g=0$ and
$\ell(I_j)=1$. Similarly, when $k=1$ then $r=d+2g-2+\ell(I_1)>1$
unless $g=0$ and $d+\ell(I_1)\le 3$. Since $\ell(I_1)\le d$ then
$\ell(I_1)=1$ and $d\le 2$.  Finally, when $k=0$ then
$r=d+2g-2$. Since there is no $d=1$ holomorphic cover of $S^2$ by a
smooth $T^2$ then $r>2$ unless $g=0$ and $d\le 2$. 

Note that since $\ell(N_j)\le \ell(I_j)$ then $n\le 2$ in all above
cases. \qed
\bigskip

\non This Lemma  motivates the following:
\begin{defn} If $k+r\le 2$, the unique element of the space $\y_{d,g}
(\ma\prod_{j=1}^k  b_{I_j}(N_j) )$ described in Lemma \ref{0.dim} 
will be called a {\em trivial cover}.
\label{D.triv}
\end{defn}
\medskip

The advantage of working with  the space $\z_{d,\chi}$ is that 
after `marking' the location of all the branch points in the target 
(which in particular means rescaling any time two of them come close 
to each other) the limit  map has no
ghost components at all and the double points of the domain occur only
above the double points of the target. This is because whenever we
start with a sequence of {\em smooth} maps there cannot be any ghost
components forming or double points appearing unless some branch
points ran into each other in the target. Therefore in this case the
limit can be thought as an {\em admissible cover} of an element of
$\ov\M_{0,k+r}$ (as described for example on p180-186 of \cite{hm}):

\begin{defn} Assume $k+r\ge 3$. The compactification 
$\;\ov\z_{d,\chi}\l(\,\ma\prod_{j=1}^k b_{I_j}(N_j)\,\r)\;$ of the space 
(\ref{def.Z}) consists of stable maps $\;\;f:C\ra A$ such that:
\begin{enumerate}
\item[(i)] the domain $C$ is a possibly disconnected curve with
 Euler characteristic $\chi$ and marked points $x_1,\dots, x_n$ so 
$st(C)\in \wt\M_{\chi,n}$;
\item[(ii)] the target $A\in\ov\M_{0,k+r}$ is a stable genus $0$ curve with
marked points $p_1,\dots,p_{k+r}$;
\item[(iii)] over the smooth part of $A$ the curve $C$ is smooth and $f$ is a
degree $d$ cover which has
ramification pattern $b_{I_i}(N_i)$ over $p_i$ for $1\le i\le k$, is simply
branched over the rest of $p_i$, $k+1\le i\le k+r$ and has no other branch
points;
\item[(iv)] the inverse image of each node of $A$ consists of nodes of $C$
with matching ramification patterns. More precisely, if  $A_1$, $A_2$ are
the two components of $A$ joined at the node $q_1=q_2$ let $C_i=f^{-1}(A_i)$
and $f^{-1}(A_1\ma\cup_{q_1=q_2} A_2)=
C_1\ma\cup_{y_i^1=y_i^2\atop i=1\dots \ell}C_2$. Then the multiplicity
$s_i$ of $f_1=f|_{C_1}$ at $y_i^1$ equals that of $f_2=f|_{C_2}$ at $y_i^2$.
\end{enumerate}
\smallskip
 
\non Let $\;\ov\y_{d,g}\l(\,\ma\prod_{j=1}^k b_{I_j}(N_j)\,\r)\;\subset
\;\ov\z_{d,\chi}\l(\,\ma\prod_{j=1}^k b_{I_j}(N_j)\,\r)\;$ denote the
corresponding compactification of the space of connected covers 
(\ref{def.Y}). 
\label{def.oy}
\end{defn} 
An element $f$ of $\oz_{d,\chi}$ is an equivalence class of triples 
consisting of the (marked) domain and target plus the covering map. 
Thus (\ref{def.diag}) extends to
\begin{equation}
\begin{array}{ccccc}
L_{p_j}&&\L_{p_j}\hskip.8in\L_{x_i}&&L_{x_i}
\\ 
\downarrow&&\searrow\hskip.6in \swarrow&&\downarrow
\\
\ov\M_{0,k+r}&\ma\longleftarrow^q&\oz_{d,\chi} 
\l(\ma\prod_{j=1}^k  b_{I_j}(N_j)\r)&@>{st}>>&
\wt\M_{\chi,n}
\end{array}
\label{def.diag2}
\end{equation}
where $\L_{x_i}\ra \oz_{d,\chi}$ and  $\L_{p_j}\ra
\oz_{d,\chi}$ are the relative cotangent bundles to the (unstabilized) domain
$C$ at $x_i$ and respectively to the target $A$ at $p_j$. Note that in
the setup above $q^*L_{p_j}=\L_{p_j}$ but in general 
$st^*L_{x_i}\ne\L_{x_i}$. This is because $A$ is a stable curve, but
$C$ might have unstable components, which get collapsed under the
stabilization map.

\bigskip

Moreover, the compactification $\oy_{d,g}$ has a natural
stratification which comes from the standard stratification of 
$\ov\M_{0,k+r}$ combined with data of the covering map which includes
the ramification multiplicity at each node of $C$ and the degree of
$f$ on each component of $C$. Each (open) stratum of the 
compactification is a smooth orbifold of (complex) dimension 
$\dim \y_{d,g}-\#\{\mbox{double points of $A$}\}$. We will show below 
that the space  $\oy_{d,g}\l(\ma\prod_{j=1}^k
b_{I_j}(N_j)\r)$ (as well as its cousin $\oz_{d,\chi}$) carries a 
fundamental class (over $\Q$) of dimension
$\max(k+r-3,0)$, which we will call a {\em ramification class}. In 
particular, the image under the stabilization map 
$st:\oy_{d,g}\l(\ma\prod_{j=1}^k b_{I_j}(N_j)\r)\ra \ov\M_{g,n}$
defines a cycle 
\best
st_*\oy_{d,g}\l(\prod_{j=1}^k b_{I_j}(N_j)\r)
\eest
on $\ov\M_{g,n}$ called a {\em ramification cycle}. We can think of 
this cycle as a 
condition on a curve $C\in\ov\M_{g,n}$, in which case it will be called a 
{\em ramification constraint}. Note that if for some $j$ we have 
$\deg(I_j)-\ell(I_j)+\ell(N_j)=0$ then the corresponding ramification 
cycle vanishes in $\ov\M_{g,n}$ by dimensional reasons.
\medskip

Moreover, let $M_j\subset N_j$ for all $j=1,\dots,k$, 
$M=\ma\sqcup_{j=1}^k M_j$  and let $\rho$, $\pi$ denote the
projections that forget those marked points which are not in $M$:
\bear
\begin{CD}
\oy_{d,g}\l(\ma\prod_{j=1}^k b_{I_j}(N_j)\r)@>{\rho}>>
\oy_{d,g}\l(\ma\prod_{j=1}^k b_{I_j}(M_j)\r)\\
\downarrow{st_n}@.
\downarrow{st_m}\\
\ov\M_{g,n}@>{\pi}>>\ov\M_{g,m}
\end{CD}
\label{big.proj}
\eear
where $m=\ell(M)$. Then $\rho$ is a finite covering map so 
 the image under $\pi_*$ of a ramification cycle in
$\ov\M_{g,n}$ is a multiple of a ramification cycle in
$\ov\M_{g,m}$. 
\medskip

Given a space $\oz_{d,\chi}\l(\ma\prod_{j=1}^k  b_{I_j}(N_j)\r)$ we can 
decompose each cover into connected components. 
In particular, for each connected component of the cover we can 
forget the marking of those points $p_j$, $j=k+1, \dots, k+r$ of 
the target over which that component is unramified. This defines a map 
$u$ which fits in the diagram 
\bear
\begin{CD}
\oz_{d,\chi}\l(\ma\prod_{j=1}^k b_{I_j}(N_j)\r)@>{u}>>
\l.\ma\bigsqcup_m \l( \bigsqcup\ma\prod_{a=1}^m \oy_{d_a,g_a}
\l( \ma\prod_{j=1}^k  b_{I_{j,a}}(N_{j,a}) \r) 
 \r)\r/S_m
\\
\downarrow{st}@.
\downarrow{\prod st}\\
\wt\M_{\chi,n}@>{=}>> \l.\ma \bigsqcup_m \l(\bigsqcup\ma\prod_{a=1}^m 
\ov\M_{g_a,n_a} \r)\r/S_m
\end{CD}
\label{z=cupy.dia}
\eear
where in the upper right hand side of the diagram the second 
union is over all (a) 
degrees $d_a\ge 1$ with $\ma\sum_{a=1}^m d_a=d$; 
(b) genera $g_a$ with $\ma\sum_{a=1}^m (2-2g_a)=\chi$; 
(c) partitions $(I_{j,a})_{a=1}^m$ of $I_j$ for each $j=1,\dots, k$;  
(d)  partitions $(N_{j,a})_{a=1}^m$ of $N_j$ for each $j=1,\dots,
k$ and (e) all possible distribution of the $r$ simple branch points 
on the connected components. As before, the 
symmetric group $S_m$ acts by permuting the $m$ domain
components. 
\bigskip

We will be mostly interested in those ramification cycles with
complicated ramification patterns only over two points. 
\begin{defn} When $k=2$ the cycle
\best
st_*\ov\y_{d,g}  \l(b_{I_1}(N_1) b_{I_2}(N_2)\r)
\eest 
on $\ov\M_{g,n}$ is called a  {\em 2-point ramification cycle}. 
\label{D.2-pt}
\end{defn}
For a 2-point ramification cycle 
$st_*\oy_{d,g}\l(b_{I_1}(N_1)b_{I_2}(N_2)\r)$ relation
(\ref{def.r}) becomes 
\best
r=2g-2+\ell(I_1)+\ell(I_2).
\eest 
So $r=0$ only for a trivial cover (see Definition \ref{D.triv}). 
For a non-trivial cover, $r\ge 1$
and 
\bear
\dim\; st_*\oy_{d,g}\l(b_{I_1}(N_1)b_{I_2}(N_2)\r) = 
2g-3+\ell(I_1)+\ell(I_2)=r-1
\label{dim.2pt}
\eear
If moreover $2g+n\ge 3$ then the codimension of 
$st_*\oy_{d,g}\l(b_{I_1}(N_1)b_{I_2}(N_2)\r)$ in $\ov\M_{g,n}$ 
(which equals the dimension of the constraint it imposes) is 
\bear
{\rm codim}\; st_*\oy_{d,g}\l(b_{I_1}(N_1)b_{I_2}(N_2)\r)= 
g+n-\ell(I_1)-\ell(I_2)=g-\ma\sum_{j=1}^2 (\ell(I_j)-\ell(N_j))
\label{codim.2pt}
\eear
In particular, in genus 0
\best
st_*\oy_{d,0}\l(b_{I_1}(N_1)b_{I_2}(N_2)\r)=0 \qquad\mbox{ if } 
\qquad \ell(I_1)+\ell(I_2)>n\ge 3. 
\eest 
More generally, when  $2g+n\ge 3$  relation (\ref{codim.2pt}) combined
with the inequalities $\ell(I_j)\ge \ell(N_j)$ and $\ell(I_j)\ge 1$
implies that 
\bear
{\rm codim}\; st_*\oy_{d,g}\l(b_{I_1}(N_1)b_{I_2}(N_2)\r) \le \min(g, g+n-2)
\label{codim.2pt.2}
\eear
For a trivial cover 
\bear
st_*\oy_{d,0}\l(b_{d}(N_1)b_{d}(N_2)\r)={1\over
d}[\ov\M_{0,n}]\in H_0(\ov\M_{0,n})\cong \Q.
\label{triv.cycle}
\eear
This follows from the diagram
\best
\begin{CD}
\oy_{d,0}(b_{d}(x_1)b_{d}(x_2)b_{1^d}(x_3))@>{\rho}>>
\oy_{d,0}(b_{d}(N_1)b_{d}(N_2))\\
\downarrow{st_1}@.
\downarrow{st}\\
\ov\M_{0,3}@>{=}>>\ov\M_{0,n}
\end{CD}
\label{big.proj.0}
\eest
after noting that the maps $\rho$ and $st_1$ have degrees $d$ and 
1 respectively.

\begin{rem} Consider the diagram (\ref{z=cupy.dia}) when  $k=2$, and
fix both a topological type for the domains of the covers in 
the moduli space 
$\oz_{d,\chi}=\oz_{d,\chi} \l( b_{I_1}(N_1) b_{I_2}(N_2)\r)$ as 
well as a particular distribution of the degree and of the branching
constraints on each component of the domain. This data picks up a 
certain component ${\cal C}$ of the moduli space $\oz_{d,\chi}$ which
is mapped by $u$ to a quotient by the symmetric group of one of the
components $\ma\prod_{a=1}^m \oy_{d_a,g_a}$. 
As usual, let $r$ be the number (\ref{def.r}) of simple branch 
points for $\oz_{d,\chi}$ and suppose that, on ${\cal C}$, $r_a$ of 
them land on the component of the cover which lies in 
$\oy_{d_a,g_a}$. In particular, $r=\ma\sum_{a=1}^m r_a$. 
But $\oz_{d,\chi}$ has dimension $\max(r-1,
0)$ while the dimension of $\ma\prod_{a=1}^m \oy_{d_a,g_a}$ is only 
$\ma\sum_{a=1}^m \max (r_a-1,0)$. Diagram (\ref{z=cupy.dia}) 
then implies that $st_*({\cal C})=0$ unless the covers in ${\cal C}$  are 
trivial on all but at most one of their connected components 
(see Definition \ref{D.triv}). Moreover, if ${\cal C}$ is a component of  
$\oz_{d,\chi}$ where all but at most one of the connected components of each  
cover are trivial, then the restriction of the map $u$ to ${\cal C}$
is an isomorphism. Therefore, the cycle $st_*\oz_{d,\chi}$ is a linear 
combination of products of 2-point ramification cycles; in each
product, all but at most one of the factors comes from a trivial
cover (see equation (\ref{triv.cycle}) for the contribution of a
trivial cover). 
\label{R.push.0} 
\end{rem} 
\bigskip

Next we describe in more detail how the strata of $\oz_{d,\chi}$ fit
together. We start with the set-theoretical picture. First notice that
there is another (coarser) stratification of $\oz_{d,\chi}$ that records 
a stratification of $\ov\M_{0,k+r}$ together with the ramification
pattern $b_S$ over the nodes of $A$ and the Euler characteristics of the
preimages of the components of $A$. Take for
example an (open) stratum where $A$ has only 2 components $A_1$ and $A_2$,
joined at the double point $q_1=q_2$. Assume moreover that the first 
$k_1$ of the points $p_i$ are on $A_1$,  the next $k_2=k-k_1$ on
$A_2$, while the remaining $r$ simple branch points are distributed in
all possible ways on the two components. Denote the closure of this 
stratum in $\ov\M_{0,k+r}$ by $D_\Gamma$ where
$\Gamma$ is the dual graph which has 2 vertices $A_i$ joined by
an edge corresponding to the node $q_1=q_2$ and tails (half edges) $p_1,\dots
,p_{k_1}$ on $A_1$ and $p_{k_1+1},\dots,p_k$ on $A_2$; sometimes we
denote this stratum by $(p_1,\dots,p_{k_1}\;|\;p_{k_1+1},\dots,p_k)$. 
 \medskip

Using the notation from  Definition \ref{def.oy}, given 
 $f\in \oz_{d,\chi}$ we start by choosing an ordering of the $\ell$ 
double points of $C$ that lie above the  node $q_1=q_2$.  We then get an 
ordered sequence $S$ of multiplicities, two smooth curves 
$C_1,\,C_2$ and two stable maps $f_i=f|_{C_i}$, $f_i:C_i\ra A_i$ such that 
\begin{enumerate}
\item[(a)] the curve $C_i$ is in $\wt \M_{\chi_i,n_i+\ell(S)}$,  
where its last $\ell(S)$ marked points are $y_1^i,\dots,y_\ell^i$; 
\item[(b)] $C=C_1\ma\cup_{y_i^1=y_i^2\atop i=1\dots \ell}C_2$ so in
particular $\chi=\chi_1+\chi_2-2\ell(S)$ and $n=n_1+n_2$;
\item[(c)] $f_1\in \z_{d,\chi_1}\l(\ma\prod_{j=1}^{k_1}
b_{I_j}(N_j)\; b_S(M^1)\r)\;$ and $f_2\in \z_{d,\chi_2}\l(b_S(M^2)\; 
\ma\prod_{j=k_1+1}^{k} b_{I_j}(N_j)\r)\;$ where $M^i=(y_1^i,\dots y_\ell^i)$.
\end{enumerate}
Consider the attaching map that (pairwise) identifies the last 
$\ell(S)$ points of $C_1$ and $C_2$
\best
\xi:\wt \M_{\chi_1,n_1+\ell(S)}\ti \wt \M_{\chi_2,n_2+\ell(S)}&\ra& 
\wt \M_{\chi,n}
\eest
given by $(C_1,C_2)\mapsto C_1\ma\cup_{y_i^1=y_i^2\atop i=1\dots
\ell}C_2$. Then all together, the data above gives a parameterization
$F$ of a stratum of $\ov\z_{d,\chi}$. More precisely, $F$ fits in the diagram
\bear
\begin{CD}
\oz_{d,\chi_1}\l(\ma\prod_{j=1}^{k_1} b_{I_j}(N_j)\; b_S\r)\;
\ma \ti\; \oz_{d,\chi_2}\l(b_S\; \ma\prod_{j=k_1+1}^{k}
b_{I_j}(N_j)\r)@>{F}>> \oz_{d,\chi}
\\
\downarrow{st\ti st}@.
\downarrow{st}\\
\wt \M_{\chi_1,n_1+\ell(S)}\ti \wt \M_{\chi_2,n_2+\ell(S)}
@>{\xi}>>\wt \M_{\chi,n}
\end{CD}
\label{def.stratum} 
\eear
where to define  $F$ we used the attaching map $\xi$
to identify the corresponding points in the inverse image of $f_1$ and
$f_2$ over $b_S$ (and thus also their images $q_1$ and $q_2$). 
As before, these points over $b_S$ are considered
marked and ordered, even though they do not appear in the
notation. The parameterization $F$ is a local embedding, but not 
necessarily injective, as the ordering of the $\ell=\ell(S)$ double
points of $C$  is not part of the original data. To keep notation
simple, we will denote by 
\bear
\oz_{d,\chi_1}\l(\ma\prod_{j=1}^{k_1} b_{I_j}(N_j)\; b_S\r)\;
\ma \ti_\xi\; \oz_{d,\chi_2}\l(b_S\; \ma\prod_{j=k_1+1}^{k}
b_{I_j}(N_j)\r)
\label{ti.xi}
\eear
the pushforward by $F$ of the fundamental class of the domain of the
parameterization (\ref{def.stratum}).  
  
\bigskip 

The inverse image 
of the stratum $D_\Gamma=(p_1,\dots,p_{k_1}\;|\;p_{k_1+1},\dots,p_k)$ of 
$\ov\M_{0,k+r}$ under $q$ can then be parameterized by
\bear
F:\bigsqcup_{\chi_i, S}\;\; 
 \oz_{\chi_1,d}\l(\ma\prod_{j=1}^{k_1} b_{I_j}(N_j)\; b_S\r)
\;\ma\ti\; \oz_{\chi_2,d}\l(b_S\; 
\ma\prod_{j=k_1+1}^{k} b_{I_j}(N_j)\r)\longra q^{-1}(D_\Gamma)
\label{def.str2} 
\eear 
where the union is over all $\chi_1$, $\chi_2$,  ordered
sequences $S$ of degree $d$  with 
$\chi=\chi_1+\chi_2-2\ell(S)$ and all possible distributions of the
$r$ simple branch points. As the target of a sequence of stable maps in
$\z_{d,\chi}$ degenerates into an element of $D_\Gamma$, the limit is
an element of $q^{-1}(D_\Gamma)$. Going backwards, we next need to understand
all possible smoothings of elements of $q^{-1}(D_\Gamma)$ into elements of
$\z_{d,\chi}$.

\bigskip

Recall that an element of $\oz_{d,\chi}$ is a triple
consisting of domain, target and a covering map. We start by looking at
smoothings of the domain and of the target. In the setup above, the normal
direction to $D_\Gamma$ inside $\ov\M_{0,k+r}$ is parameterized by the
line bundle $ \L_{q_1}^*\otimes \L_{q_2}^*$ whose fiber at 
$A_1\ma \cup_{q_1=q_2} A_2$ is   $T_{q_1}A_1\otimes
T_{q_2}A_2$.  Similarly, the normal bundle of the $\ell$-nodal stratum in 
$\wt \M_{\chi,n}$ is 
$\ma\bigoplus_{i=1}^{\ell}\; \L_{y_i^1}^*\otimes \L_{y_i^2}^*
$, whose fiber at $C=C_1\ma\cup_{y_i^1=y_i^2\atop i=1\dots \ell}C_2$ is
$\ma\bigoplus_{i=1}^{\ell} \;T_{y^1_i}C_1\otimes T_{y^2_i}C_2$. 
However, for a fixed smoothing $A_\la$ of $A$ not all smoothings
$C_\mu$ of $C$ give rise to a stable map $f:C_\mu\ra A_\la$; here 
$\la\in T_{q_1}A_1\otimes T_{q_2}A_2$ and $\mu=(\mu_1,\dots \mu_\ell)$
with  $\mu_i\in  T_{y_i^1}C_1\otimes T_{y_i^2}C_2$. This can be best
seen in local coordinates $z_{m,i}$ at $y_m^i$ and $w_i$ at $q_i$. In
these coordinates 
\bear
w_i=f_i(z_{m,i})=a_{m,i} \cdot (z_i)^{s_m}+\mbox{ higher order}
\label{def.a}
\eear 
while the smoothings of $C$ and $A$ are given by
$z_{m,1}\cdot z_{m,2}=\mu_m$, $m=1,\dots, \ell$ and $w_1\cdot
w_2=\la$. Therefore $f:C\ra A$ can be extended to a smooth cover
$f_{\mu,\la}:C_\mu\ra A_\la$ only when 
\best
\la =a_{m,1}\,a_{m,2}\,\mu_m^{s_m}\qquad \mbox{ for all }m=1,\dots, \ell
\eest
to highest order. For example this fact is proven (in a more general setting)
using PDE methods in \cite{ip4}. It was also stated in the original 
Harris-Mumford paper \cite{hmu}. Moreover, in the algebraic-geometrical
setting, the deformation argument of Caporaso  and Harris  \cite{ch}
could be extended to this case. After all, in \cite{ch} they have
studied stable maps into $\P^2$ with prescribed  contact constraints 
along a line $L$, and the case above is simply  a dimensional
reduction where the pair $(\P^2,L)$ gets replaced by $(\P^1,p)$.  

\bigskip

Summarizing, given a pair
$(f_1,f_2)$ in the domain of the parameterization  (\ref{def.stratum}), 
equation (\ref{def.a}) defines a canonical section
\bear
\si_q:\z_{d,\chi_1}(\dots b_S)\ma\ti \z_{d,\chi_2}(b_S\dots) 
\longrightarrow
\ma\bigoplus_{i=1}^\ell\;\l( \L_{x_{n_i}}\otimes 
\L_{y_{n_i}}\r)^{s_i}  \otimes 
\l(\L^*_{q_1}\otimes \L^*_{q_2}\r)
\label{lead.term.2}
\eear
given by $\si_q=(a_1,\dots, a_\ell)$ with $a_m=a_m^1\cdot a_m^2$. For a
fixed smoothing of the target $\la\in  \L^*_{q_1}\otimes \L^*_{q_2}$ the
possible smoothings of the domain $\mu=(\mu_1,\dots, \mu_\ell)$
correspond to solutions of the equations 
\bear
\la=a_1\mu_1^{s_1}=\dots =a_\ell \mu_\ell^{s_\ell} 
\label{la=amu}
\eear
There are $|S|=\prod s_i$ many such solutions, differing by roots of
unity.  This describes the local model in the normal direction to a 
stratum parameterized by 
$\z_{d,\chi_1}(\dots b_S)\ma\ti \z_{d,\chi_2}(b_S\dots)$
inside the compactification $\oz_{d,\chi}$. Moreover, this shows that
as cycles, the pullback  of $D_\Gamma$  is
\bear
q^*(D_\Gamma)=\bigsqcup_{\chi_i, S} \;\;{|S|\over \ell(S)!}\;\; 
\oz_{\chi_1,d}\l(\ma\prod_{j=1}^{k_1} b_{I_j}(N_j)\; b_S\r)
\;\ma\ti_\xi\; 
\oz_{\chi_2,d}\l(b_S\; \ma\prod_{j=k_1+1}^{k} b_{I_j}(N_j)\r)
\label{def.str3} 
\eear 
where the union is over all  $\chi_1$, $\chi_2$, ordered sequences
$S$ of degree $d$ with $\chi=\chi_1+\chi_2-2\ell(S)$ and all possible 
distributions of the $r$ simple branch points. The $1\over \ell(S)!$
weight comes from the fact that the ordering of the $\ell(S)$ double
points of $C$ is not part of the original data of an element in 
$\ov \z_{d,\chi}$. 

\begin{rem} Note that the solution space to the equations
(\ref{la=amu}) has several branches intersecting at the origin 
(which corresponds to the boundary stratum) so the compactification 
$\oz_{d,\chi}$ described in  Definition \ref{def.oy} is not in general
an orbifold. However, it can be desingularized by including as part of
the data besides the triple $f:C\ra A$ a choice of roots of unity for
the leading term section (\ref{lead.term.2}). This desingularized 
compactification becomes then a version of the
space of {\em twisted covers} defined in \cite{av}. In any event, we
will only use the fact that (each component of) $\oz_{d,\chi}$ carries a 
fundamental class (with rational coefficients) and so $st_*[\oz_{d,\chi}]$
defines a class on $\wt\M_{\chi,n}$.  
\end{rem}
\medskip

As a particular case of (\ref{def.str3}) we get the following 
\begin{theorem} Let $q:\oz_{d,\chi}\l(\ma\prod_{j=1}^{k} b_{I_j}(N_j)\r)\ra
\ov\M_{0,k+r}$ be as in (\ref{def.diag2}) and let $D_\Gamma$ be the
codimension one stratum of $\ov\M_{0,k+r}$ where the first $k_1$
points are on a bubble, the next
$k_2=k-k_1$ points are on a different bubble and the remaining $r$ points
are distributed all possible ways. Then as cycles in $\ov \M_{\chi,n}$ 
\bear
st_*q^*(D_\Gamma)=\sum  \;{|S|\over \ell(S)!}\;\; 
st_*\l(\; \oz_{\chi_1,d}\l(\ma\prod_{j=1}^{k_1} b_{I_j}(N_j)\; b_S\r)
\;\ma\ti_{\xi}\; \oz_{\chi_2,d}\l(b_S\; \ma\prod_{j=k_1+1}^{k} 
b_{I_j}(N_j)\r) \;\r)
\label{deg.form}
\eear
where the sum is over all $\chi_1$, $\chi_2$, ordered sequences $S$ 
of degree $d$ with $\chi=\chi_1+\chi_2-2\ell(S)$ and all possible 
distributions of the $r$ simple branch points.

When $q$ is restricted to the space $\oy_{d,g}$ of connected covers 
then we get cycles in $\ov\M_{g,n}$ and in the sum above we
keep only those configurations of domains $C_1, \;C_2$ whose image under
the attaching map $\xi$ is connected.   
\label{T.main.form}
\end{theorem}

\smallskip

Note that the  equal sign in (\ref{deg.form}) is only an equality in
homology, because the proof in \cite{ip5} (which is done in the
symplectic category) only shows that the compactification 
$\oz_{\chi,n}$ is diffeomorphic to the local model (\ref{la=amu}). 
However, an algebraic-geometrical proof of the local model
(\ref{la=amu}) would give the equality in the Chow ring. 

\bigskip

\begin{ex} Suppose $k=2$, $k_1=k_2=1$ and $2g+n\ge 3$. The right hand side of
(\ref{deg.form}), when restricted to connected genus $g$ covers, 
involves terms of type
\bear
\oz_{\chi_1,d}\l( b_{I_1}(N_1)\; b_S\r)
\;\ma\ti_{\xi}\; \oz_{\chi_2,d}\l(b_S\; b_{I_2}(N_2)\r) 
\label{term.0.1}
\eear
The pushforward by $st$ of such term, using relation (\ref{ti.xi}) 
and diagram (\ref{def.stratum}),   is equal to 
\best
\xi_*\l(st_*\oz_{\chi_1,d}\l( b_{I_1}(N_1)\; b_S\r)
\;\ti\; st_* \oz_{\chi_2,d}\l(b_S\; b_{I_2}(N_2)\r)  \r)
\eest
By Remark \ref{R.push.0}, each component of $st_*\oz_{\chi_i,d}$ 
is a multiple of a product of 2-point ramification cycles; the 
factors in the product correspond to (unstabilized) domain
components. Moreover, the discussion following 
Definition \ref{D.2-pt} implies that on all genus 0 components the 2-point
ramification cycles either vanish or else are multiples of the
fundamental class. Suppose we fix a topological type of the
(unstabilized) domain, and a fixed distribution of the ramification
patterns on each component of the domain. Then the pushforward by $st$
of  the corresponding component of (\ref{term.0.1}) equals a rational 
multiple of products of 2-point ramification cycles on the components 
of the stabilized domain.  More precisely,  suppose the
stabilized domain consists of components of genus $g_a$ with $n_a$ 
special points labeled by $M_a$, all glued together according to the 
dual graph via the attaching map
\bear
\xi:\prod_{a=1}^h \ov\M_{g_a,n_a} \ra \ov\M_{g,n}
\label{attach.xi}
\eear
By convention, when the dual graph has no edges (i.e. stabilized
domain is smooth), the attaching map is the identity. Then the
component of the right hand side of (\ref{deg.form}) 
corresponding to the attaching map (\ref{attach.xi}) is a 
linear combination (with rational coefficients) of terms of type 
\bear
\xi_*\l( \prod_{a=1}^h st_{*}\oy_{d_a,g_a}\l( b_{I_{a1}}(N_{a1})\;
b_{I_{a2}}(N_{a2}) \r)\r)\in H_*(\ov\M_{g,n})
\label{prod.st}
\eear
where $N_{a1}\sqcup N_{a2}=M_{a}$. Note that the term (\ref{prod.st}) 
vanishes unless on all genus 0 components $\ell(N_{a1})=\ell(I_{a1})$ and 
$\ell(N_{a2})=\ell(I_{a2})$ (see relation 
(\ref{codim.2pt})). Moreover, since all terms in (\ref{deg.form}) are
codimension one, then only terms of type (\ref{prod.st}) 
for which the domain has at most one node can appear (with nonzero
coefficient) in the right hand side of (\ref{deg.form}). 
\end{ex}

\begin{defn} Consider ramification cycles
$
{\cal C}_a=st_*\oy_{d_a,g_a}\l(\ma\prod_{i=1}^{k_a} b_{I_{ai}}(N_{ai})\r)
$ 
on $\ov\M_{g_a,n_a}$  where $2g_a+n_a\ge 3$. For each attaching map
$\xi$ as in (\ref{attach.xi}) the cycle
$
\xi_* ( \ma\prod_{a=1}^h  {\cal C}_a) 
$
is called a {\em generalized ramification cycle} on
$\ov\M_{g,n}$. In particular, such a cycle for which $k_a=2$ for all
$a=1,\dots, h$ will be  called a  generalized 2-point ramification cycle.
\end{defn}
With this definition, Theorem \ref{T.main.form} implies in particular
that when $D_\Gamma$ is a codimension one boundary stratum of
$\ov\M_{0,r+k}$, then $st_*q^*D_\Gamma$ is a linear combination of 
codimension one generalized 2-point ramification cycles. 
\begin{defn} Let $\Theta$ be a linear combination of generalized 
ramification cycles on  $\ov\M_{g,n}$. Those terms of $\Theta$ which 
are constructed using the attaching map of a boundary stratum of 
$\ov\M_{g,n}$ will be called {\em lower order terms}. The sum of the 
other terms forms the {\em symbol of} $\Theta$. By 
convention, if all the terms are lower order, we take the symbol to be 0. 
\end{defn} 

\begin{rem}  For $M\subset \{x_1\dots,x_n\}$ consider the map 
$\pi:\ov\M_{g,n}\ra \ov\M_{g,m}$ that forgets the marked points which
are not in $M$. Suppose ${\cal C}$ is a generalized ramification cycle 
on $\ov\M_{g,n}$. Since the attaching maps commute with the forgetful
maps, diagram (\ref{big.proj}) implies that $\pi_*{\cal C}$ is a 
(rational) multiple of a  generalized ramification cycle on 
$\ov\M_{g,m}$. Note that even if  ${\cal C}$ is nonzero,  $\pi_*{\cal
C}$ might vanish (by dimensional reasons for example). 
\label{R.forget}
\end{rem}

\begin{rem} Theorem \ref{deg.form} generalizes to higher codimensional
boundary strata in $\ov\M_{0,k+r}$. In particular, let $D$ denote the 
codimension $m-1$ boundary strata of $\ov\M_{0,2+r}$ consisting of 
linear chains of $m$ $\P^1$'s such that $p_1$ is on the 
first bubble, $p_2$ on the last 
bubble and the other $r$ points $p_j$, $j=3,\dots, r+2$ are
distributed in some fixed way on the $m$ components. Then 
\best
q^{-1}(D)\subset 
\oz_{d,\chi}(b_{I_1}(N_1)b_{I_2}(N_2))
\eest
is similarly parameterized by a disjoint union of spaces
\bear
\oz_{d,\chi_1}(b_{I_1}(N_1) b_{S_1})\ma\ti_{\xi_1} \oz_{d,\chi_2}(b_{S_1}
b_{S_2})\ma\ti_{\xi_2}\dots \ma\ti_{\xi_{m-1}}  
\oz_{d,\chi_m}(b_{S_{m-1}} b_{I_2}(N_2))
\label{gen.z}
\eear
where each attaching map $\xi_i$ identifies the corresponding points
over $b_{S_i}$ for $i=1,\dots, m-1$. So $st_*q^*(D)$ can also be
written as a linear combination of generalized 2-point ramification 
cycles of codimension $m-1$. 
\label{R.gen.form}
\end{rem}

\medskip

Next, fix a moduli space $\oy_{d,g}\l(\ma\prod_{j=1}^k
b_{I_j}(N_j)\r)$ such that $2g+n\ge 3$ (so in particular $k+r\ge 3$ by
Lemma \ref{0.dim}). 
Assume in what follows that the point  $x_i$ has 
prescribed ramification index $s_i$ and image $p_j$. We can
consider  the `universal family'
$\oy_{d,g}\l(b_{1^d}(x_{0})\ma\prod_{j=1}^k b_{I_j}(N_j)\r)$ obtained
by adding extra marked points $x_{0}$ to the domain and $p_0$ to the
target, together with the diagram
\bear
\begin{CD}
@.\ov\M_{g,n+1}@<{st_{n+1}}<<
\oy_{d,g}\l(b_{1^d}(x_{0})\ma\prod_{j=1}^k b_{I_j}(N_j)\r)@>{q_{n+1}}>>
\ov\M_{0,k+r+1}\\
@.x_i\uparrow\downarrow{\pi_1}@.\downarrow{\pi_0}@.{p_j}\uparrow
\downarrow{\pi_2}\\
@.\ov\M_{g,n}@<{st_n}<<\oy_{d,g}\l(\ma\prod_{j=1}^k b_{I_j}(N_j)\r)@>{q_n}>>
\ov\M_{0,k+r}
\end{CD}
\label{big.dia}
\eear
where $\pi_0$ is the map that forgets both the marked point $x_{0}$ 
on the domain and its image $p_{0}$ on the target.  The images of the
canonical sections $x_i$, $p_j$ are the strata 
$D_{0,i}\subset \ov\M_{g,n+1}$ and respectively $D_{0,j}\subset 
\ov\M_{0,r+1}$ where $x_0$ and $x_i$ and respectively $p_0$ and $p_j$
are the only marked points on a genus 0 bubble. 
\medskip

The covers in the preimage $q_{n+1}^{-1}(D_{0,j})\subset 
\oy_{d,g}\l(b_{1^d}(x_{0})\ma\prod_{j=1}^k b_{I_j}(N_j)\r)$ have a
very special form. Because on the genus 0 bubble containing $p_0$ and
$p_j$ there are no other branch points, then over this component the 
cover consists of $\ell(I_j)$ spheres totally ramified above $p_j$ and 
$p_\infty$ ($p_\infty$ is the double point of the target where the 
bubble is attached). Only the sphere that contains $x_0$ is nontrivial, 
the rest are trivial covers (see Definition \ref{D.triv}). 
When the point $x_0$ is on the
same bubble as the point $x_i$ we denote the corresponding canonical 
section by 
\bear
\si_i:\oy_{d,g}\l(\prod_{j=1}^k b_{I_j}(N_j)\r)\ra 
\oy_{d,g}\l(b_{1^d}(x_{0})\prod_{j=1}^k b_{I_j}(N_j)\r)
\eear
and let $\Si_i$ denote its image. Then 
$st_{n+1}\circ \si_i=x_i\circ st_n$ and $q_{n+1}\circ  \si_i=p_j\circ
q_n$ with the notations of (\ref{big.dia}). In particular, 
this discussion shows that
\bear
\pi_1^* st_{n*}\oy_{d,g}\l(\ma\prod_{j=1}^k b_{I_j}(N_j)\r)= 
st_{n+1*}\oy_{d,g}\l(b_{1^d}(x_{0})\ma\prod_{j=1}^k b_{I_j}(N_j)\r)
\label{pull.back.0}
\eear 
Moreover, 
\begin{lemma}  Consider the space 
$\oy_{d,g}=\oy_{d,g}\l(\ma\prod_{j=1}^k b_{I_j}(N_j)\r)$ where all the 
preimages of  all the marked points of the target 
(including all branch points) are marked.  Suppose moreover that $x_i$ 
is a marked point in the domain with image $p_j$ and ramification
index $s_i$. If $L_{x_i}\ra \ov\M_{g,n}$ and $L_{p_j}\ra
\ov\M_{0,r}$  are the  relative cotangent bundles to the domain and
respectively the target then  over $\oy_{d,g}$ we have
\bear
st^* L_{x_i}^{s_i}\;= \;q^* L_{p_j}
\label{eq.RH}
\eear
\end{lemma}
\pf Since  $L_{x_i}=x_i^*\O(-D_{0,i})$ and 
$L_{p_j}=p_j^*\O(-D_{0,i})$ then
\best
st_n^*L_{x_i}=st_n^* x_i^*\O(-D_{0,i})= \si_i^*st_{n+1}^*\O(-D_{0,i})\\
q_n^*L_{p_j}=q_n^* p_j^*\O(-D_{0,j}) = \si_i^*q^*_{n+1}\O(-D_{0,j}).
\eest
But all the points over $p_j$ are marked so all the covers
in $q_{n+1}^{-1}(D_{0,j})$ have domains with  $x_0$ and at least one of the
other points over $p_j$ on the same bubble. Moreover, the only instance
where $x_0$ and $x_i$ are the only two marked points on a genus 0
bubble are those covers in $\Si_i$. Then  (\ref{deg.form}) implies
that
\best
q^*_{n+1}\O(-D_{0,j})=\O(-s_i\Si_i) \qquad \mbox{ along } \Si_i
\eest
where $s_i$ is the ramification index of point $x_i$. The condition 
that all the preimages of all the marked points of
the target (including all branch points) are marked implies in
particular that all the domains of the covers are stable curves and
therefore 
\best
st_{n+1}^*\O(-D_{0,i})=\O(-\Si_i) \qquad \mbox{ along } \Si_i.
\eest
Combining  the last four displayed equations we then get (\ref{eq.RH}).\qed

\bigskip
\section{Polynomials in descendants}
\medskip

In this section we describe how to express a product of 
$\psi_i=c_1(L_{x_i})$  classes on $\ov\M_{g,n}$ 
(or more precisely the intersection product of their Poincare duals) 
as a linear combination of generalized ramification cycles. 
\medskip

The basic idea is simple: to begin with we choose a 2-point ramification 
cycle $\oy_{d,g}$ so that the map 
$st:\oy_{d,g}\ra \ov\M_{g,n}$  is of finite (nonzero)
degree. Then we use equation (\ref{eq.RH}) to relate 
$st^*L_{x_i}\ra \oy_{d,g}$
to the pull back $q^*L_{p_j}$ of the relative cotangent bundle 
$L_{p_j}$ to the target $\P^1$ at $p_j$, the image of $x_i$ under 
the covering map.  But we know that the Poincare dual of 
$c_1(L_{p_j})$ is a codimension 1 boundary cycle $D_\Gamma$ in 
$\ov\M_{0,r}$. Then Theorem \ref{T.main.form} implies that the
Poincare dual of $\psi_i$ is linear combination of generalized 
2-point ramification cycles on $\ov\M_{g,n}$. 

\bigskip
In what follows the descendant on the target $c_1(L_{p_j})$ will 
be denoted  by $\wt \psi_j$ to avoid confusing it with the descendant
on the domain $\psi_j=c_1(L_{x_j})$. Also, in the rest of the paper, 
we will often add or forget marked points. Note to begin with that if  
$\pi_{0}:\ov \M_{g,n+1}\ra \ov \M_{g,n}$ 
is the map that forgets the marked point $x_0$ then
\bear
\psi_i=\pi_{0}^*\psi_i +D_{i,0}
\label{phi.pull.back}
\eear
where $D_{i,0}$ is the boundary strata in $\ov \M_{g,n+1}$ 
consisting of domains where $x_i$ and $x_0$ are the {\em only} points 
on a $g=0$ 
bubble.  Similarly, for tautological classes we have
\bear
\kappa_i=\pi_{0}^* \kappa_i + \psi_{0}^i
\label{kappa.pull.back}
\eear
Moreover, if $\xi$ is the attaching map (\ref{attach.xi}) of a
boundary stratum of $\ov\M_{g,n}$ then the pullback by $\xi$ of the
relative cotangent bundle to $x_i$ is the relative cotangent bundle to
$x_i$, so 
\bear
\xi^*\psi_1=\psi_1.
\label{xi.psi}
\eear

\begin{ex} Let us illustrate the procedure described at the beginning of this
section on the following example: when $g=1$ it is known that 
$\psi_1=\delta_0/12$ in $\ov\M_{1,1}$, where $\delta_0$ is the
boundary stratum which corresponds to
a nodal sphere. To see this, we start by  writing  any element in 
$\ov\M_{1,1}$ as a degree 2 cover of $\P^1$ branched at 4 points such 
that the marked point $x_1$ is one of the branch points. As long as
the branch points are not ordered, such cover is in fact unique.  
Fix now 2 other branch points and let 
$\oy_{2,1}\ma=^{\rm def}\oy_{2,1}(b_2(x_1)b_2^2)$ denote the
corresponding space of covers, and let  
\best
\oy_{2,1}(b_2(\psi_1)b_2^2)\ma=^{\rm def} st^*\psi_1\cap [\oy_{2,1}]
\eest 
denote the Poincare dual of $st^*\psi_1$ in $\oy_{2,1}$. 
Then the stabilization map $st:\oy_{2,1}\ra \ov \M_{1,1}$ is a
degree $3\cdot 2=6$ cover (for each possible choice of the 2 out of 3 
remaining branch points) so
\bear
st_*\oy_{2,1}(b_2(\psi_1)b_2^2)=6 \psi_1
\label{g=1.1}
\eear
Next, relation (\ref{eq.RH}) gives
$2 st^*\psi_1=q^*\wt\psi_1$ on $\oy_{2,1}$. On the other hand, on
$\ov\M_{0,4}$ $\wt\psi_1$ is Poincare dual to the boundary stratum 
$D_\Gamma$ which consists of $p_1,\;p_2$ on one bubble and 
$p_3,\;p_4$ on the other so
\bear
2\oy_{2,1}(b_2(\psi_1)b_2b_2)= q^*(D_\Gamma)
\label{g=1.2}
\eear 
Now use the degeneration formula (\ref{deg.form}). Since the degree
is 2 and total genus is 1, the only term that can appear is $S=(1,1)$
with genus 0 on both sides, i.e. 
\best
st_*q^*(D_\Gamma)=st_*\l( \oy_{2,0}(\,b_2(x_1)b_2b_{1,1}\,)\ma \ti_\xi
\oy_{2,0}(\,b_{1,1}b_2^2\,)\r)
\eest
But there is only one genus 0 degree two map, and under the
stabilization map the component on the right gets collapsed to a
point, and therefore
\bear
st_*q^*(D_\Gamma)=\delta_0
\label{g=1.3}
\eear
Combining (\ref{g=1.1}), (\ref{g=1.2}) and (\ref{g=1.3}) gives the
relation $\psi_1=\delta_0/12$.
\end{ex}  
\smallskip

More generally,
\begin{theorem} Assume $g\ge 1$, $n\ge 1$ and $n+g\ge 3$. 
Then the Poincare dual of any degree $m$ monomial in descendant 
classes on $\ov\M_{g,n}$ can
be written as a linear combination of generalized 2-point ramification
cycles on $\ov\M_{g,n}$, coming from a cover of degree at most
$d=g+n-1$. The nonzero terms appearing in the symbol are codimension
$m$ cycles of type 
\best
st_*\oy_{a,g} (b_{I_1}(N_1)b_{I_2}(N_2)) 
\eest 
where $a\le d$, $N_1\sqcup N_2=\{x_1,\dots, x_n\}$ and
$\ell(I_1)+\ell(I_2)=g+n-m$.  
\label{T.1}
\end{theorem}

Note that $\ell(N_j)\le \ell(I_j)$ so adding we
get $n\le g+n-m$. In particular, the Theorem implies that when 
$m\ge g+1$  or $m\ge g+n-1$ there are no nonzero terms in the symbol,
and so the degree $m$  monomial in descendant classes vanishes when
restricted to $\M_{g,n}$. 

Moreover, a closer analysis of the proof of Theorem \ref{T.1} shows
that the terms appearing in the symbol have either $\ell(I_1)=1$ or 
$\ell(I_2)=1$. But since this is irrelevant for this paper, 
we leave the details to the reader. 
\medskip

\non{\bf Proof of Theorem \ref{T.1}.} 
Consider the ramification cycle (as defined in Section 1)
\best
\oy_{d,g,n}=\oy_{d,g}(b_{1^d}(N)b_{d})
\eest
where $N=(x_1,x_2,\dots, x_n)$. Under the assumptions of the Theorem, 
when $d=g+n-1$ Lemma \ref{al.ne.0} below shows that 
$st:\oy_{d,g,n}\ra \ov \M_{g,n}$
is map of finite, {\em nonzero} degree $\deg(st)\ne 0$.
\medskip

Now let $\psi_1^{m_1}\dots \psi_n^{m_n}$ be a monomial on 
$\ov\M_{g,n}$ of degree $m=\sum m_j\ge 0$. The Poincare dual of 
$st^*\l(\psi_1^{m_1}\dots \psi_n^{m_n}\r)$ in 
$\oy_{d,g}(b_{1^d}(N)\,b_{d}\,)$ is  
\best
st^*\l( \psi_1^{m_1}\dots\psi_n^{m_n}\r) 
\; \cap\; [\oy_{d,g}(b_{1^d}(N)\,b_{d}\,)] 
\eest
so the Poincare dual 
\bear
PD\l(\psi_1^{m_1}\dots\psi_n^{m_n}\r)\; =\; (\deg(st))^{-1} \cdot 
st_*\l(   st^*\l( \psi_1^{m_1}\dots\psi_n^{m_n}\r) 
\; \cap\; [\oy_{d,g}(b_{1^d}(N)\,b_{d}\,)]\r) 
\label{phi.1}
\eear
The Theorem then follows by induction on the degree $m$ of the monomial 
$\psi_1^{m_1}\dots \psi_n^{m_n}$. The case $m=0$ comes directly 
from relation (\ref{phi.1}). Now suppose the result is true for $m-1$, 
so we need to prove it for $m$. Consider a monomial $\psi_1^{m_1}\dots
\psi_n^{m_n}$ of degree $m\ge 1$. Without loss of generality we may
assume that $m_1\ge 1$. Then relation (\ref{phi.1}) implies  
\bear
PD\l(\psi_1^{m_1}\dots\psi_n^{m_n}\r)
\; =\; (\deg(st))^{-1} \;\psi_1 \cap  
st_* \l(st^*( \psi_1^{m_1-1}\dots\psi_n^{m_n}) \cap 
[\oy_{d,g}(b_{1^d}(N)\,b_{d}\,)]\r)
\label{psi.2}
\eear
By induction, $st_* \l(st^*( \psi_1^{m_1-1}\dots\psi_n^{m_n}) \cap 
[\oy_{d,g}(b_{1^d}(N)\,b_{d}\,)]\r)=(\deg st)\cdot  
\psi_1^{m_1-1}\dots\psi_n^{m_n}$  is a
linear combination of generalized 2-point ramification cycles. Thus 
the cycle (\ref{psi.2}) is a linear
combination of terms of type
\best
\psi_1\cap  \xi_{*}\l( \prod_{a=1}^m st_* \oy_{d_a,g_a}\l
( b_{I_{a1}}(N_{a1}) b_{I_{a2}}(N_{a2})\r) \r)= 
\xi_{*}\l( \xi^*\psi_1\cap \prod_{a=1}^m st_* \oy_{d_a,g_a}
\l( b_{I_{a1}}(N_{a1}) b_{I_{a2}}(N_{a2})\r) \r)
\eest
Using relation (\ref{xi.psi}) and applying  Lemma
\ref{L.InductiveStep} to the factor containing the 
marked point $x_1$ then completes the inductive step. 
\qed

\begin{lemma} Let  $d=g+n-1$. Then the degree of the map 
\best
st: \oy_{d,g}(b_{1^d}(N)b_{d}) \ra \ov\M_{g,n}
\eest 
is nonzero as long as $g\ge 1$, $n\ge 1$ and $g+n\ge 3$. Moreover, the
degree of $st$ vanishes when $g=0$ or $n=0$ or $g=n=1$.
\label{al.ne.0}
\end{lemma}  
\pf We begin by noting that when $d=g+n-1$, dimension count shows
that the domain and target of the map $st$ have the same
dimension. The vanishing part of the Lemma follows immediately after 
noting that when $g=0$
the domain of $st$ is empty (since $\ell(N)>d$), while when $n=0$ or 
$g=n=1$ the fiber of $st$ is one dimensional. 

For $d=g\ge 2$ Mumford proved in Section 7 of \cite{mu} that
the degree of the stabilization map $st:\oy_{d,g}(b_{d}) \ra \ov\M_{g,0}$ is
nonzero. In particular, this implies that the degree of the map 
$st:\oy_{d,g}(b_d b_{1^d}(x_1))\ra \ov\M_{g,1}$ is nonzero as well, 
because once we write a generic Riemann surface as an element of 
$\oy_{d,g}(b_d)$, adding a generic marked point $x_1$ gives an element of 
 $\oy_{d,g}(b_d b_{1^d}(x_1))$. This
proves the Lemma in the case $n=1$ and $g\ge 2$. 

The case when $n\ge 2$ and $g\ge 1$ follows by methods similar to 
those of Section
5 of \cite{hmu}. More precisely, fix a 
generic (smooth) genus $g$ Riemann surface $C$ with $n$ marked points 
$x_i$, $i=1,\dots,n$. 
It is enough to show that we can find $g$ points 
$y_0,\dots, y_{g-1}$ on $C$ such that 
\best
\sum_{i=1}^n x_i +\sum_{i=1}^{g-1}y_i \sim d y_0
\eest  
Then as long as $g\ge 1$ and $n\ge 2$, a dimension count shows that the
points $y_0,\dots,y_{g-1}$ are distinct and distinct from the
points $x_1, \dots, x_n$, thus producing the required degree $d$
cover. To show existence,  let  $J(C)$ be the 
Jacobian of $C$, $u:C\ra J(C)$ be the 
Abel-Jacobi map, and  $C_d=Sym^d(C)$. Consider the maps 
$v:C\ra J(C)$ and $w:C_{g-1}\ra J(C)$ given by $v(y)=u(dy)=du(y)$ and 
$w(D)=u(D)+u(\ma\sum_{i=1}^{n}x_i)$. We need to show that the
intersection between the image of $v$ and that of $w$ is nonempty.
But the image of $w$ is a translate of the $\Theta$ divisor and moreover
$v^*w_*[C_{g-1}]=v^*([\Theta])=dg\ne 0$.  \qed

\begin{lemma} Fix a 2-point ramification cycle 
$ st_*  \oy_{d,g}( b_I(N)b_J(M))$ on $\ov\M_{g,n}$ 
where $N\sqcup M=\{ x_1, \dots, x_n\} $. Then the cycle
\best
\psi_1\cap st_*  \oy_{d,g}( b_I(N)b_J(M))= 
st_* \l( st^*\psi_1 \cap \oy_{d,g}( b_I(N)b_J(M)) \r) 
\eest
can be written as a linear combination of generalized 2-point 
ramification cycles; its symbol consists of terms of type 
\best
st_*\oy_{a,g}( b_{I_1}(N_1)b_{J_1}(M_1))
\eest
where $a\le d$, $N_1\sqcup M_1=\{x_1,\dots, x_n\}$ and 
$\ell(I_1)+\ell(J_1)=\ell(I)+\ell(J)-1$. 
\label{L.InductiveStep}
\end{lemma}
\pf The result is trivially true when $r=0$, i.e. 
$\oy_{d,g}( b_I(N)b_J(M))$ is zero dimensional (see Lemma
\ref{0.dim}). So we may assume $r>0$.  

The first step is  to replace $st^*\psi_1$ by
a multiple of $q^*\wt\psi_1$, where $\wt \psi_1=c_1(L_{p_1})$ is the
first Chern class of the relative cotangent bundle to the target $\P^1$ at
$p_1$. For that, we temporarily mark the location of the other
$\ell(I)-\ell(N)$ points in the preimage of $p_1$, $\ell(J)-\ell(M)$
points in the preimage of $p_2$ and each of the $d-1$ points in the
preimage of each of the other $r=2g-2+n-\ell(I)-\ell(J)$ simple 
branch points. All together, we add $b=r(d-2)+2g-2$ extra marked points,
getting a corresponding 2-point cycle $\oy_{d,g,n+b}$ in which all the
preimages of all the branch points are marked. Consider the diagram
\begin{equation}
\begin{array}{cccl}
\ov \M_{g,n+b}&\ma\longleftarrow^{st_b}& \oy_{d,g,n+b}
\\
\pi_b \downarrow&&\rho_b \downarrow&\searrow^q
\\
\ov \M_{g,n}&\ma\longleftarrow^{st}& \oy_{d,g,n}&\ma\ra^{q}
\ov\M_{0,2+r}
\end{array}
\end{equation}
where $ \oy_{d,g,n}=\oy_{d,g}(b_I(N)b_J(M)$. Then 
$\rho_b: \oy_{d,g,n+b}\ra  \oy_{d,g,n}$ has 
finite {\em nonzero} degree $\deg(\rho_b)=
(\ell(I)-\ell(N))!\cdot (\ell(J)-\ell(M))!\cdot (d-1)!^r\ne 0$ so
\best
(\deg \rho_b)\cdot st^*\psi_1=\rho_{b*}\,\rho_b^* \,st^*\psi_1=
\rho_{b*}\,st_b^*\,\pi_b^*\psi_1 
\eest
Moreover, the stabilization map 
$st_b:\oy_{d,g,n+b}\ra \ov \M_{g,n+b}$ does not collapse any
components of the domain. Therefore, the relative cotangent bundle 
$\L_{x_1}\ra \oy_{d,g,n+b}$ to the domain is equal to
the pullback by $st_b$ of $L_{x_1}\ra\ov\M_{g,n+b}$. Using  
formula (\ref{phi.pull.back}) repeatedly and pulling back by $st_b$
gives then the relation
\best
st_b^*\,\pi_b^*\psi_1= c_1(\L_{x_1})- st_b^*\D_1
\eest
on  $\oy_{d,g,n+b}$, where $\D_{1}=\ma\sum_{L} D_{1,L}$ and $D_{1,L}$ 
is the boundary
strata in $\ov \M_{g,n+b}$ where the marked point $x_1$ and a subset 
$L$ of the $b$ new marked points are the only points on a $g=0$
bubble.  

\medskip
Now on  $\oy_{d,g,n+b}$ all the preimages of the marked points of the
target are marked so the relation (\ref{eq.RH})  
implies that $\L_{x_1}^{s_1}=q^*L_{p_1}$ so
\best
c_1(\L_{x_1})={1\over s_1}\cdot q^*(c_1(L_{p_1}))={1\over s_1}\cdot 
q^*(\wt\psi_1)
\eest 
Combining the last three displayed equations we get
\bear
\psi_1\cap  \oy_{d,g,n}(b_I(x_1)b_J(M))=
{1\over s_1\deg(\rho_b)}\cdot st_*q^*(\wt\psi_1)- 
{1\over \deg(\rho_b)}\cdot \rho_{b*}st_b^*\D_1
\label{rel.1}
\eear
Next, we use the fact that in
$\ov\M_{0,2+r}$ we have $r\cdot \wt\psi_1=D$ where 
$D=\ma \sum_{j=3}^{r+2}D_{\Gamma_j}$ 
and $D_{\Gamma_j} $ is the boundary strata that has the marked point
$p_1$ on a bubble and $p_2$, $p_j$ on a different bubble, while the
remaining $r-1$ branch points are distributed all possible
ways. Note that in $D$ the  strata which has a bubble
containing $p_1$ and precisely $r_1$ of the points $p_j$ with $j\ge
3$ appears with coefficient $r_2=r-r_1$. Applying the degeneration formula 
(\ref{deg.form}) for each $j$ and summing then gives
\bear
st_*q^*(\wt\psi_1)={1\over r}\; st_*q^*(D)=
\sum {|S|\over \ell(S)!}\cdot {r_2\over r}\cdot 
\;st_*\l(\oz_{d,\chi_1}( b_{I}(N) b_S)\ma\ti_\xi \oz_{d,\chi_2}(b_Sb_J(M))
\r)
\label{deg.1}
\eear
where the sum is over all $\chi_1$, $\chi_2$, $r_1$, $r_2$, 
ordered sequences $S$ such that $\deg S=d$,
$\chi_1+\chi_2-2\ell(S)=2g-2$, $r_1+r_2=r$,  
over all possible identifications that lead to a connected domain and
over all possible distributions of the $r$ simple branch points such
that $r_1$ are on the left component. In any case, this show that the
first term on the right hand side of equation (\ref{rel.1}) is a
linear combination of generalized 2-point ramification cycles. 

On the other hand $\rho_{b*}st^*_b \D_1$ is also equal to a  linear 
combination of similar generalized 2-point ramification cycles. 
This is because $st_b$ doesn't collapse any components, thus 
$st_b^*\D_1$ consists of stable maps in $\oy_{d,g,n+b}$ whose domain 
is an element of $\D_1$. In particular, the target of these maps must 
be a bubble tree with $p_1$ on one side and $p_2$ on the other. 

Using (\ref{rel.1}) and (\ref{deg.1}) we then conclude that on 
$st_*\oy_{d,g}( b_I(N)b_J(M))$, $\psi_1$ can be written as a linear 
combination of the generalized 2-point ramification classes. 
The statement about the structure of the
symbol follows immediately by a dimension count. \qed  
\bigskip

Because of Remark \ref{R.forget}, an immediate consequence of Theorem 
\ref{T.1} is the following: 
\begin{cor} Assume $g\ge 2$, $g+n\ge 3$ and let 
$\Pi_k:\ov\M_{g,n}\ra \ov \M_{g,n-k}$ be a forgetful map. Then the
Poincare dual of the class $\Pi_{k*}\l(\psi_1^{m_1}\dots \psi_n^{m_n}\r)$
 on $\ov\M_{g,n-k}$ can be written as linear combination 
of generalized 2-point ramification cycles on 
$\ov\M_{g,n-k}$ whose symbol consists of codimension $\ma\sum_{j=1}^n
m_j -k$ terms of type
\best
st_*\oy_{a,g}( b_{I_1}(N_1)b_{I_2}(N_2))
\eest
where $a\le d$, $N_1\sqcup N_2=\{x_1,\dots, x_{n-k}\}$ and 
$\ell(I_1)+\ell(I_2)= g+n-\ma\sum_{j=1}^n m_j$. 
\label{P.induction2}
\end{cor}
Note that since $\ell(N_j)\le \ell(I_j)$ then in particular the symbol 
vanishes when $\ma\sum_{j=1}^n m_j>g+k$.

\begin{rem} If one is interested not only in the shape of the symbol,
but in the actual formula then it is convenient to start with a cover of
degree as small as possible, so there would be fewer terms to
consider. In this context, one can use  the fact that any complex
structure can be written as a degree $d=\l[ {g+1\over 2}\r]+1$ cover 
of $\P^1$ to adapt the proof of Theorem \ref{T.1} to get the 
following:
\begin{prop} Any polynomial in descendant classes on $\ov\M_{g,n}$ can be 
written as a linear combination of generalized  ramification 
constraints coming from covers of degree at most $\l[ {g+1\over 2}\r]+1$.
\label{T.2}
\end{prop}
For example, when $g$ is odd, one would start with the space 
$\oy_{g,d}\l(\ma\prod_{i=1}^n b_{1^d}(x_i)\r)$ for which the degree of
the stabilization map is nonzero, while when $g$ is 
even, one would use instead the space $\oy_{g,d}\l(b_{2,1^{d-2}}(x_1)
\ma\prod_{i=2}^n b_{1^d}(x_i)\r)$. Then one uses the fact that in 
$\ov\M_{0,r+n}$ the Poincare dual of any monomial in descendant 
classes $\wt\psi_j$ can be expressed as a linear combination of 
boundary strata corresponding to linear chains of $\P^1$'s. 
In the end, after using Remark \ref{R.gen.form}, one would get generalized 
ramification cycles with at most two complicated branch points on 
each component of the target (but not technically 2-point ramification 
cycles, because of the presence of constraints of type $ b_{1^d}(x_i)$).  
\end{rem}

\medskip

\subsection{Proof of Theorem \ref{T.main}} 
\bigskip

Suppose we start with a degree $m$ monomial in $\psi$ and $\kappa$  
classes on $\ov \M_{g,n}$. Then using the formulas (\ref{phi.pull.back}) and 
(\ref{kappa.pull.back}) we can express any such polynomial 
as a linear combination of terms of type
\best
\Pi_{k*}(\psi_1^{m_1}\dots\psi_{n+k}^{m_{n+k}})
\eest
for some $k$'s, where $\Pi_{k}:\ov\M_{g,n+k}\ra \ov\M_{g,n}$ is the map that 
forgets the last $k$ marked points, and $\ma\sum_{j=1}^{n+k} m_j=m+k$. 
For example,
\best
\kappa_a\kappa_b=\Pi_{2*}(\psi_1^{a+1}\psi_2^{b+1})-\Pi_{1*}(\psi_1^{a+b+1})
\eest
It is therefore enough to prove Theorem \ref{T.main} for classes of
type
\bear
\Pi_{k*}\l(\psi_1^{m_1}\dots \psi_{n+k}^{m_{n+k}}\r)\in H^{m-k}
(\ov\M_{g,n})
\label{big.mono}
\eear
where $m=\ma\sum_{i=1}^{n+k} m_i\ge g+k$.  We actually prove that 
the Poincare dual of such class can be written as a linear combination 
of generalized ramification cycles with vanishing symbol on 
$\ov\M_{g,n}$, i.e. all terms are  coming from the boundary 
$\partial \ov\M_{g,n}$.

Corollary \ref{P.induction2}, with $n$ replaced by $n+k$, implies that 
the Poincare dual of the class (\ref{big.mono}) can be written as a 
linear combination of generalized ramification 
 cycles whose symbol consists of terms of type  
\best
\oy_{a,g}( b_{I_1}(N_1)b_{I_2}(N_2))
\eest
with $a\le d$, $N_1\sqcup N_2=\{x_1,\dots, x_n\}$ and
$\ell(I_1)+\ell(I_2)=g+n+k-m$. So when
$m\ge g+k$ we have
\best
n=\ell(N_1)+\ell(N_2)\le \ell(I_1)+\ell(I_2)=g+n+k-m\le n
\eest 
thus all terms in the symbol vanish unless  $\ell(N_j)=\ell(I_j)$ for
$j=1,2$  and $m=g+k$.
\smallskip

When $n\le 1$ there are no such terms since $\ell(I_j)\ge 1$, so the
symbol vanishes.  Moreover, note that when $n=0$ even for $m=g+k-1$ a
similar string of inequalities shows that the symbol also vanishes,
 implying Looijenga's result \cite{l1} (in homology). 
\smallskip

When $n\ge 2$,  Proposition \ref{P.B.is.bd} below completes  the 
proof of Theorem \ref{T.main}.  \qed

\begin{prop} Suppose $g\ge 1$ and $\ell(I_i)=\ell(N_i)$ for
$i=1,2$. Then the codimension $g$ cycle  on $\ov\M_{g,n}$  
\best
{\cal C}=st_* \oy_{d,g}(b_{I_1}(N_1)b_{I_2}(N_2))
\eest 
can be written as a linear combination of generalized ramification 
cycles of type 
$\xi_* (\ma\prod_{a=1}^h \pi^*_a {\cal C}_a ) $
where $\xi:\ma\prod_{a=1}^h \ov\M_{g_a,n_a}\ra \ov\M_{g,n}$ is the
attaching map of some {\em boundary strata} of $\ov\M_{g,n}$, 
$\pi_a:\ov\M_{g_a,n_a}\ra \ov\M_{g_a,m_a}$ is a forgetful map
(this includes the identity map in the case $m_a=n_a$) and 
${\cal C}_a$ is a 2-point ramification  cycle on 
$ \ov\M_{g_a,m_a}$ coming from a
degree $d_a\le d$ cover.  
In particular, the symbol of this linear combination vanishes. 
\label{P.B.is.bd}
\end{prop}
\pf  We prove the statement by induction on both the degree $d$ and
the number of marked points $n$. It is enough to prove that the cycle 
${\cal C}$ can be written as a linear combination of cycles of type 
$\xi_* (\ma\prod_{a=1}^h \pi^*_a {\cal C}_a) $ which either come 
from the boundary or else have only one component (i.e. $h=1$) and
for this component either $d_1<d$ or $m_1<n_1=n$. 
\smallskip

Assume $x_1\in N_1$ and let $N_1'=N_1\setminus \{x_1\}$. 
Consider the cycle 
\best
\oy_{d,g}(B_{I_1}(N_1')B_{1^d}(x_1)B_{2,1^{d-2}}B_{I_2}(N_2))
\eest
which corresponds to fixing the location of the marked points
$p_1,\dots, p_4$ on the target. But on $\ov\M_{0,4+r}$ the divisor
corresponding to fixing the location of $p_1,\dots,p_4$ is linearly 
equivalent to the boundary stratum $D=(p_1p_2|p_3p_4)$ where $p_1, p_2$
are on a bubble and $p_3, p_4$ are on a different bubble. For
simplicity we denote $q^* (p_1p_2|p_3p_4)=
\oy_{d,g}( b_{I_1}(N_1')b_{1^d}(x_1)\;|\;b_{2,1^{d-2}}b_{I_2}(N_2))$. 
Since the 
stratum $(p_1p_2|p_3p_4)$ is linearly equivalent to the stratum 
$(p_1p_3|p_2p_4)$ then
\bear
st_*\oy_{d,g}(b_{I_1}(N_1') b_{1^d}(x_1)\;|\;b_{2,
1^{d-2}}b_{I_2}(N_2))=
st_*\oy_{d,g}( b_{I_1}(N_1') b_{2,1^{d-2} }\;|\;b_{1^d}(x_1)b_{I_2}(N_2))
\label{4.pt}
\eear
as codimension $g$ cycles in $\ov\M_{g,n}$. But the degeneration formula 
(\ref{deg.form}) implies that both sides of  (\ref{4.pt}) are
linear combination of pushforwards by $st$ of terms of type
\bear
&&\oz_{d,\chi_1}
(b_{I_1}(N_1')b_{1^d}(x_1)b_S)\ma\ti_\xi\oz_{d,\chi_2}(b_Sb_{2,
1^{d-2}}b_{I_2}(N_2))
\quad\mbox{ and }
\label{term.1}
\\
&&\oz_{d,\chi_1}( b_{I_1}(N_1')b_{2, 1^{d-2}}b_S)
\ma\ti_\xi\oz_{d,\chi_2}(b_Sb_{1^d}(x_1)b_{I_2}(N_2))
\label{term.2}
\eear
respectively. We need to show that the only term not lying in the
boundary of $\ov\M_{g,n}$ and with $d_1=d$, $m_1=n_1$ is the term
${\cal C}$; moreover ${\cal C}$ should appear in (\ref{4.pt}) with
nonzero coefficient. Let ${\cal C}'$ be such a term appearing after
stabilization in the symbol of (\ref{term.1}) or (\ref{term.2}).  This
means that before stabilization we have a degree $d$ genus $g$
component on one side and all the components on the other side are
genus $0$ totally ramified over the node of the target; otherwise
collapsing them would produce a double point of the (stabilized)
domain. Moreover, before stabilization we can have at most one marked
point on each genus 0 component (since when $g\ge 1$ the strata of
$\ov\M_{g,n}$ having stable $g=0$ components are in the boundary).

\smallskip

Suppose first that ${\cal C}'$ appears in the symbol of (\ref{term.1}). We
have two cases to consider:
\begin{enumerate}
\item[(a)] the genus $g$ component is on the left. But since
$\ell(I_2)=\ell(N_2)$ the genus 0 component on the right which
contains the simple ramification point cannot be totally ramified over
$p_4$ so will have to contain two of the marked points in $N_2$,
contradiction.
   
\item[(b)] the genus $g$ component is on the right. Since 
$\ell(I_1)=\ell(N_1')+1$ there can be at most one genus 0
component which is not totally ramified over $p_1$ (otherwise two of
the points in $N_1'$ would be on the same genus 0 component). But one
of the genus 0 components must also contain $x_1$, so the only
possibility is if all genus 0 components were totally ramified over
$p_1$ and moreover $x_1$ would be on the only genus 0 component not
containing a point from $N_1'$. After pushing forward by $st_*$ this
term contributes
\best
s_1 st_* (b_{I_1}(N_1)b_{I_2}(N_2)) =s_1 {\cal C}
\eest 
to the right hand side of (\ref{4.pt}), where $s_1$ is the multiplicity 
of $x_1$ in $I_1$. 
\end{enumerate}

\non Next suppose that ${\cal C}'$ appears in the symbol of 
(\ref{term.2}). We also have two cases to consider:
\begin{enumerate}
\item[(a)] the genus $g$ component is on the left. Since
$\ell(I_2)=\ell(N_2)$ then each genus 0 component on the right has at
least one of the marked points of $N_2$. But one of these  genus 0
components must also have $x_1$, contradiction.  

\item[(b)] the genus $g$ component is on the right. Since 
$\ell(I_1)=\ell(N_1')+1$ there can be at most one genus 0
component which is not totally ramified over $p_1$, and this component
can have at most 2 points over $p_1$ (otherwise two of the points in 
$N_1'$ would land on the same genus 0 component). This genus 0
component must contain the simple ramification point and only one of
the points $x_a\in N_1'$, the other point over $p_1$ being
unmarked. The order of ramification over the node of the target of this
component must then be equal to the sum of the multiplicities of the
points over $p_1$. Denote by $\wh I$ the sequence obtained from
$I_1$ by erasing the multiplicity  corresponding to $x_1$
and adding it to the multiplicity corresponding to $x_a$. After 
collapsing the genus 0 components this term is equal to a multiple of 
\best
st_*\oy_{d,g}\l( b_{\wh I}(N_1') b_{1^d}(x_1) b_{I_2}(N_2) \r)
\eest
where $\ell(\wh I)=\ell(N_1')$. By relation (\ref{pull.back.0}) this term is
equal to $
\pi_1^* st_*\oy_{d,g}\l( b_{\wh I}(N_1') b_{I_2}(N_2) \r)$ where 
 $\pi_1$ is the map that forgets the marked point $x_1$. Therefore 
it is pulled back from a moduli space with fewer marked points. 
\end{enumerate}
This concludes the inductive step and with it the proof of Proposition
\ref{P.B.is.bd}.\qed
\bigskip

We finish this paper by proving the following result, which was 
recently conjectured by  Vakil \cite{v} for the Chow group. 
\begin{prop} The Poincare dual of any degree $m$ monomial in 
descendant or tautological classes on $\ov\M_{g,n}$ can be written as
a linear combination of classes coming from the stratum of 
$\ov\M_{g,n}$ which has at least $m+1-g$ genus 0 components.  
\label{P.ravi}
\end{prop}
\pf  The result is already known in genus 0 or 1, so we prove it for $g\ge
2$. As in the proof of Theorem \ref{T.main}, Corollary \ref{P.induction2}
implies that the Poincare dual of 
any degree $m$ monomial in $\psi$ and $\kappa$ classes can be written as a 
linear combination of codimension $m$ generalized 2-point ramification 
cycles. Each such generalized 2-point ramification cycle is of type 
$\xi_*( \ma\prod_{a=1}^m {\cal C}_a )$ where 
$\xi:\ma\prod_{a=1}^m\ov\M_{g_a,m_a}\ra \ov\M_{g,n}$ is the
attaching map of some stratum of $\ov\M_{g,n}$ (including possibly the
top stratum) and each ${\cal C}_a$ is a 2-point ramification cycle of type 
${\cal C}_a=st_* \oy_{d_a,g_a}(b_{I_{a1}}(N_{a1})b_{I_{a2}}(N_{a2}))$. 
The codimension of such ${\cal C}_a$  is at most $g_a$ by  relation
(\ref{codim.2pt.2}). But by induction (on the dimension of the moduli 
space $\ov\M_{g,n}$) we can prove that any 2-point ramification cycle 
${\cal C}=st_* \oy_{d,g}(b_{I_{1}}(N_{1})b_{I_{2}}(N_{2}))$ 
can be written as a linear combination  of 
generalized 2-point ramification cycles of type  
\bear
\xi_*\l( \ma\prod_{a=1}^m {\cal C}_a \r),\qquad \mbox{ where  } \qquad 
{\cal C}_a=  \pi_a^* st_* \oy_{d_a,g_a}(b_{I_{a1}}(N_{a1})b_{I_{a2}}(N_{a2}))
\label{xi*}
\eear
where moreover ${\rm codim }\; {\cal C}_a\le g_a-1$ on all 
$g_a\ge 1$ components. This is because either ${\cal C}$ already has
codimension less then $g$ or else Proposition \ref{P.B.is.bd} shows
that it can be written as a linear combination of generalized ramification
cycles of type (\ref{xi*}) coming from a boundary strata (in which
case each ${\cal C}_a$ comes from a lower dimensional moduli space). 

Therefore the Poincare dual of any degree $m$ monomial in 
$\kappa$ and $\psi$ classes can be written as a linear combination 
of codimension $m$ generalized ramification cycles of type (\ref{xi*})
for which ${\rm codim }\; {\cal C}_a\le g_a-1$ on all $g_a\ge 1$ components.  
Fix such a codimension $m$ generalized ramification cycle. We only
need to show that the domain of the corresponding attaching 
map $\xi$ has at least $m+1-g$ genus 0 components.  
Let $k$ be the number of double  points and  
$\ell$ be the number of genus 0 components of the corresponding 
stratum of $\ov\M_{g,n}$. Then 
\best
m=k+\sum_{a=1}^m \mbox{ codim }{\cal C}_a\le  
k+\sum_{g_a\ge 1} (g_a-1)= k+\sum_{a=1}^m(g_a-1)+\ell=g-1+\ell
\eest
where the last equality follows from the Euler characteristic relation
$2-2g=\ma\sum_{a=1}^m (2-2g_a) -2k$. Therefore
$\ell \ge m+1-g$. \qed
\bigskip

\bigskip



\end{document}